\documentclass[11pt,final,3p,times,authoryear]{elsarticle}

\usepackage{amssymb}

\usepackage{algorithm2e}
\usepackage{algorithmic}
\RestyleAlgo{ruled}
\SetKwComment{Comment}{}{ }

\usepackage{amsmath}

\usepackage{geometry}
\usepackage[symbol]{footmisc}
\usepackage{amsthm}
\usepackage{bm}
\usepackage{amsfonts}
\usepackage{amssymb}
\usepackage{color}
\usepackage{graphicx}
\graphicspath{{}}
\usepackage{natbib}
\usepackage{float}
\usepackage{wrapfig}
\usepackage[hidelinks]{hyperref}
\usepackage{siunitx}
\usepackage{caption}
\usepackage[dvipsnames]{xcolor}
\usepackage{subcaption}
\usepackage{mathtools}
\usepackage[english]{babel}
\usepackage{booktabs}
\usepackage{longtable}
\allowdisplaybreaks
\usepackage[flushleft]{threeparttable}
\usepackage{soul}
\usepackage{makecell}
\usepackage{setspace}
\usepackage{enumitem}
\newcommand{\midsize}{\fontsize{11pt}{11pt}\selectfont}

\newcommand\blfootnote[1]{%
  \begingroup
  \renewcommand\thefootnote{}\footnote{#1}%
  \addtocounter{footnote}{-1}%
  \endgroup
}

\begin{document}

\begin{frontmatter}

\title{\LARGE{Two-Echelon Delivery Vehicle Sharing and Repositioning in Hyperconnected Urban Logistic Networks}}

\author[label1]{\Large{Xiaoyue Liu} \footnote[1]{Corresponding author. Address: 755 Ferst Drive NW, Atlanta, GA, 30332, United States. Tel: +1 470-819-6831.}}
\author[label1,label2]{Benoit Montreuil\blfootnote{Email addresses: xliu800@gatech.edu (X. Liu), benoit.montreuil@isye.gatech.edu (B. Montreuil).}}

\affiliation[label1]{organization={Physical Internet Center, H.Milton Stewart School for Industrial and Systems Engineering},
            addressline={Georgia Institute of Technology}, 
            city={Atlanta},
            postcode={30332}, 
            state={GA},
            country={United States}}
\affiliation[label2]{organization={Coca-Cola Material Handling \& Distribution Chair},
            country={United States}}

\begin{abstract}
In response to the growing demand for sustainable and efficient urban deliveries, this study introduces a two-echelon vehicle sharing and repositioning problem for containerized delivery operations within a hyperconnected urban logistics system. We leverage a Physical Internet (PI)-enabled three-tier logistic hub network, comprising gateway, local, and access hubs, to facilitate efficient flows. By adopting containerized delivery, vehicles can rapidly swap standardized modular containers at hubs to reduce handling time. Moreover, inspired by the PI concept of open resource sharing, we determine optimal service routes within a two-echelon structure that jointly utilizes heterogeneous vehicle fleets and enables dynamic vehicle relocation across hubs. We formulate this problem as a multi-period integer program that integrates path-based service vehicle planning with arc-based container routing. To address real-world large-scale instances, we propose a decomposition-based heuristic with capacity-aware flow assignment, which partitions the problem into subproblems structured by echelons and regions. A case study on the Atlanta metropolitan area demonstrates the effectiveness of the proposed model and solution approach. Experimental results show that the two-echelon hyperconnected delivery system reduces CO$_2$-eq emissions by 45.0\% and total costs by 16.8\% at full market share compared to a traditional single-echelon alternative, while enabling vehicle repositioning further lowers costs by up to 17.7\%.
\vspace{0.2in}
\end{abstract}


\end{frontmatter}

\newpage
\section{Introduction}
\label{sec:intro}
\subsection{Motivation}
\label{sec:mtv}

Over the past decade, the rapid expansion of e‑commerce has driven unprecedented growth in the urban delivery industry. In 2024, U.S. e‑commerce sales reached $1.19$ trillion, a new peak that accounted for 16.1\% of total retail transactions \citep{SellersCommerceUSEcommerce}. This surge in demand, coupled with high population densities in metropolitan areas, has intensified urban freight traffic, leading to worsening emissions and congestion that undermine urban livability \citep{schoder2016impact, gonzalez2023exploring}. According to the \cite{WEF2024UrbanLogistics}, the number of delivery vehicles in cities could increase by 61\% by 2030, driving a 60\% rise in delivery-related carbon emissions, and a 14\% increase in congestion, with commuters losing up to five additional minutes per trip. Ensuring the sustainability of urban logistics is no longer just desirable but has become a critical necessity. Despite this mounting pressure, traditional urban logistics systems remain largely inefficient and fragmented. Static routing, isolated operations, and lack of coordination between carriers often result in overlapping service areas, low vehicle utilization, and excessive vehicle mileage \citep{Caldarone2024LastMile}. These systemic inefficiencies, combined with rising consumer expectations for speed and flexibility \citep{Shankar2025LastMileDeliveryCosts}, have pushed last-mile delivery costs to over 53\% of total shipping expenses \citep{WEF2024UrbanLogistics}. Therefore, realizing efficient and sustainable urban logistics requires a more fundamental rethinking of network design and operational strategies.

To address these challenges, the Physical Internet (PI) initiative, proposed by \cite{montreuil2011toward}, offers a promising solution. Drawing inspiration from the Digital Internet's data transmission, PI envisions an open and interconnected logistics system enabling seamless movement and efficient consolidation of physical goods through standardized, modular containers traversing shared networks and infrastructures \citep{montreuil2013foundations}. Key concepts of PI include the use of PI containers, protocol-based routing, optimized asset utilization, and open resource sharing across hyperconnected logistic hub networks, aiming to achieve improvements in efficiency and sustainability \citep{wu2025towards}. Building on this foundation, \cite{crainic2016physical} introduced the hyperconnected city logistics, outlining nine key principles as a PI framework for designing sustainable urban logistics systems. Subsequently, \cite{crainic2020planning} developed an optimization model for tactical planning in a two-tier hyperconnected network, providing a method to evaluate the benefits of the framework. Several recent works have explored the network design \citep{ muthukrishnan2021potential, kulkarni2022resilient}, shared capacity deployment strategies \citep{faugere2022dynamic, liu2023logistics, liu2025dynamic}, resource planning and allocation methods \citep{li2024stochastic, liu2024dynamic, xu2024network}, and container flow routing planning \citep{kim2021hyperconnected, liu2025network}, demonstrating promising reductions in cost and carbon emissions across PI-enabled hyperconnected logistics systems.

\begin{figure}[h!]
    \centering
    \includegraphics[width=16.cm]{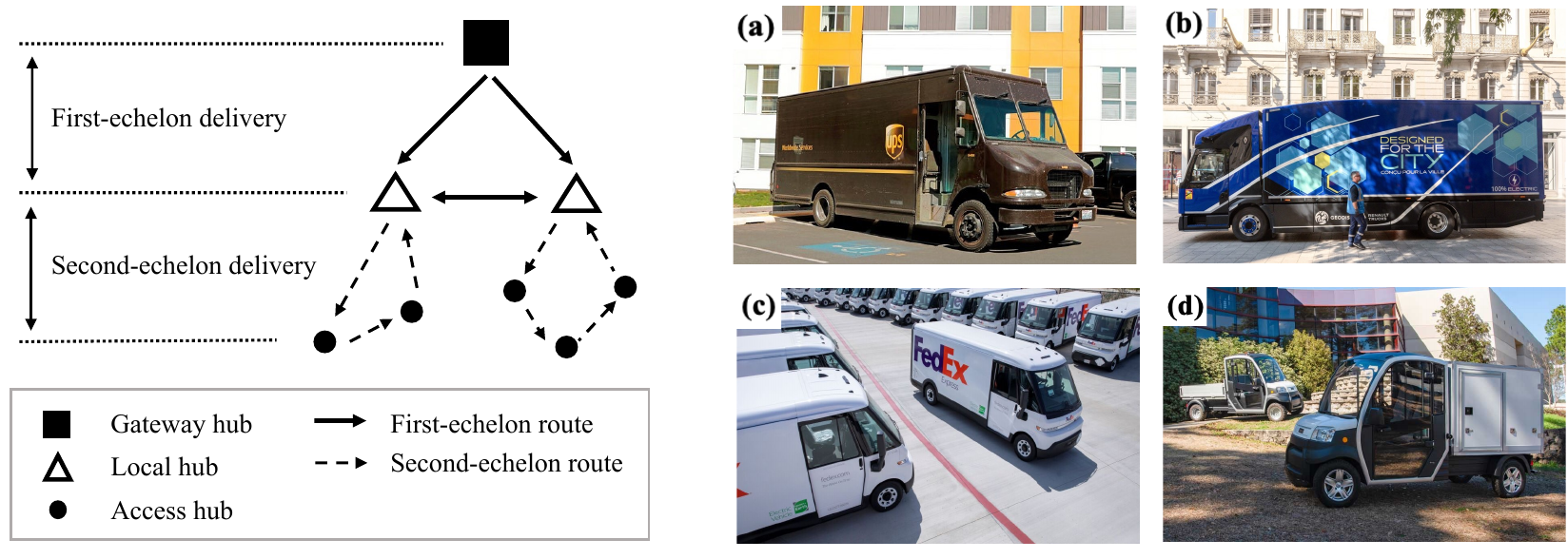}
    \caption{{Illustration and examples of two-echelon urban delivery vehicles: (a) \textit{UPS Trucks} (www.ups.com), (b) \textit{GEODIS} (geodis.com), (c) \textit{FedEx} (www.fedex.com), and (d) \textit{Club Car} (www.clubcar.com).}}
    \label{fig:two_echelon_vehicles}
    \vspace{-0.1in}
\end{figure}

Inspired by the open resource sharing and PI-container concepts from the Physical Internet, this paper introduces the dynamic vehicle sharing and repositioning problem for containerized urban deliveries. This new operational problem aims to determine both vehicle-service planning and container flow routing over multiple planning periods, serving customer demand in urban delivery contexts. Accordingly, this problem presents a multi-period setting and belongs to the NP-hard class as it inherits the complexity of the vehicle routing problem \citep{cordeau2007vehicle}. In this study, we utilize a three-tier hyperconnected urban logistic network composed of gateway, local, and access hubs to integrate and support freight movements, aligning with the network structure described in \citet{montreuil2018urban}. In a real-world context, gateway hubs represent large distribution centers on the urban fringe near major highways for regional freight consolidation, local hubs are smaller intra-city facilities for sorting and transferring shipments to specific districts, and access hubs serve as hyper-local endpoints, such as automated parcel lockers, situated within neighborhoods for convenient customer pickup. Building upon this network, we consider two‑echelon containerized deliveries where standardized containers are swapped directly between vehicles at local hubs, thereby eliminating handling delays and facilitating flow transfer. Unlike classical vehicle routing models, we allow for dynamic vehicle sharing and repositioning among same‑level hubs in response to spatial-temporal demand fluctuations. Correspondingly, the service capacity of the system is dynamically managed through decisions on the number of scheduled service routes activated and the allocation of available delivery vehicles across the hub network. Moreover, in order to adapt to city streets, our two-echelon routing strategy employs the coordinated use of vehicles in two primary sizes. As illustrated in Figure \ref{fig:two_echelon_vehicles}, the first-echelon (FE) utilizes high-capacity delivery trucks (e.g., (a) and (b)) for longer-distance transfers between gateway and local hubs; while the second-echelon (SE) employs smaller delivery vans (e.g., (c) and (d)) for agile, localized distribution between local and access hubs. 

\vspace{0.1in}

\subsection{Contribution}
\label{sec:ctb}

Overall, the key contributions of the paper are summarized as follows:
\renewcommand\labelitemi{\tiny$\bullet$} 
\begin{itemize}[itemsep=0.25mm,parsep=0.1mm]
    \item We introduce a new operational planning problem: the two-echelon vehicle sharing and repositioning problem tailored for sustainable containerized deliveries in hyperconnected urban logistic networks.
    \item We develop a multi-period hybrid optimization model that combines path-based service vehicle planning with arc-based container flow routing over a three-tier logistic hub network.
    \item We propose an effective solution approach leveraging an echelon-region decomposition heuristic with a tailored rolling horizon framework and capacity-aware flow assignment to solve large-scale freight delivery problems in megacities.
    \item We validate the proposed framework through a real-world case study using freight flow data from the Atlanta metropolitan area. The results demonstrate the benefits of the vehicle sharing and repositioning strategy, two-echelon hyperconnected delivery system, and alternative-fuel technologies in enhancing cost efficiency, operational performance, and environmental sustainability.
\end{itemize}

The rest of this paper is structured as follows. Section \ref{sec:lr} provides a literature review of related work. In Section \ref{sec:pdf}, we define the problem and formulate an arc- and path-based hybrid optimization model. Section \ref{sec:sa} details the proposed solution approach. Section \ref{sec:cs} presents our case study and derives insights from numerical results. Finally, Section \ref{sec:con} concludes the paper and discusses future research directions.

\section{Literature review}
\label{sec:lr}
The Vehicle Routing Problem (VRP), first introduced by \cite{dantzig1959truck}, has long been recognized as a fundamental problem in the field of logistics and transportation research. It focuses on determining optimal vehicle routes to serve a set of dispersed customers while minimizing total cost or distance. However, the classic VRP, which typically assumes direct deliveries from a central depot, faces significant challenges in the context of modern urban logistics \citep{perboli2011two}. The increasing congestion in city areas, route restrictions for heavy-duty trucks, and negative environmental impacts associated with long-distance last-mile deliveries render the single-echelon approach inefficient and unsustainable \citep{zeng2014hybrid}. These limitations have highlighted the need for more sophisticated network structures that can better address the operational complexity of urban freight distribution. To address these challenges, the Two-Echelon Vehicle Routing Problem (2EVRP) has emerged as a more realistic and efficient paradigm. The 2EVRP explicitly models a two-tier distribution structure: Goods are first transported from central depots to intermediate facilities, known as satellites, using large-capacity vehicles (referred to as first echelon); From these satellites, smaller and more agile vehicles perform the final deliveries to end customers (referred to as second echelon). This hierarchical structure enables better zone-based deliveries and low-emission urban mobility \citep{heidari2023green}. Comprehensive surveys, including \citet{cuda2015survey} and \citet{sluijk2023two}, have documented the growing academic interest in the 2EVRP and its extensions.

To better reflect real-world operational complexities, recent research has incorporated additional features into the 2EVRP framework, as summarized in Table \ref{tab:lr_problem}. One such extension is the inclusion of multi-period (MP) planning. Studies such as \cite{crainic2016modeling}, \cite{wang2021two}, \cite{mohamed2023two}, and \cite{al2024multi} have adopted this setting to model more dynamic logistics environments. For example, \cite{mohamed2023two} integrated strategic location planning with operational second-echelon routing decisions over a multi-period planning horizon. While \cite{al2024multi} studied a multi-period multi-commodity 2EVRP, which serves dynamic customer-specific demands in urban networks. Another recent development involves containerized delivery (CD), where goods are consolidated into standardized containers to facilitate transshipment between vehicles and enhance handling efficiency. \cite{muhlbauer2021parallelised, mohri2024designing} investigated van-bike delivery systems that transship standardized containers from vans to cargo-bikes at satellites, and through case studies in Munich and Melbourne, demonstrated the operational and environmental advantages of such containerized two-echelon urban logistics. \cite{achamrah2024gradient} extended the literature by studying a stochastic 2EVRP that incorporates both delivery and pickup of containers with cargo bikes. In addition, recent research has begun to explore vehicle repositioning (VR), which allows vehicles to be dynamically relocated between satellites and leverages resource sharing among logistic hubs over time. For instance, \cite{crainic2016modeling} extended the two-echelon framework by allowing vehicles to travel to external zones with surplus demand, pick up additional loads, and continue serving new zones. Similarly, \cite{yu2020two} incorporated the open vehicle routing problem into the 2EVRP, where a mix of closed company-owned vehicles and open rented vehicles is jointly deployed. However, as far as we know, there has been limited effort in integrating the three settings within the 2EVRP framework. Our study aims to fill this gap in the existing literature by developing a model that simultaneously incorporates these features in urban delivery planning.

\vspace{0.05in}
\begin{table}[h!]
\centering
\caption{{Characteristics of the related studies of the Two-Echelon Vehicle Routing Problem (2EVRP).}} \label{tab:lr_problem}
\resizebox{\textwidth}{!}{
\begin{threeparttable}
  \begin{tabular}{lcccccccccc}
    \toprule
    & \multicolumn{3}{c}{\textbf{Problem settings} \tnote{a}} & \multicolumn{2}{c}{\textbf{Method} \tnote{b}} & \multicolumn{4}{c}{\textbf{Case study} \tnote{c}} \\\cmidrule(lr){2-4} \cmidrule(lr){5-6} \cmidrule(lr){7-10}
    \textbf{Reference} & \textbf{MP} & \textbf{CD} & \textbf{VR} & \textbf{Form.} & \textbf{Sol.} & \textbf{Real} & \textbf{\#Depots} & \textbf{\#Satellites} & \textbf{\#Customers} \\ 
    \midrule
    \cite{crainic2011multi} & & & & A & H & & 1 & 5 & 50 \\ 
    \cite{perboli2011two} & & & & A & H & & 1 & 5 & 50 \\ 
    \cite{hemmelmayr2012adaptive} & & & & A & H & & 1 & 10 & 200 \\ 
    \cite{nguyen2012solving} & & & & A & H & & 1 & 10 & 200 \\ 
    \cite{baldacci2013exact} & & & & P & E & & 1 & 6 & 100 \\ 
    \cite{zeng2014hybrid} & & & & P & H & & 1 & 5 & 50 \\ 
    \cite{crainic2016modeling} & \checkmark & & \checkmark & P & H & & 2 & 3 & 25 \\ 
    \cite{breunig2019electric} & & & & P & H & & 1 & 10 & 200 \\ 
    \cite{yu2020two} & & & \checkmark & A & H & & 1 & 5 & 100 \\ 
    \cite{wang2021two} & \checkmark & & & A & H & \checkmark & 1 & 18 & 300 \\ 
    \cite{muhlbauer2021parallelised} & & \checkmark & & P & H & \checkmark & 1 & 15 & 300 \\ 
    \cite{dellaert2021multi} & & & & H & E & & 3 & 5 & 100 \\ 
    \cite{lv2022two} & & & & A & H & & 1 & 25 & 519 \\ 
    \cite{mohamed2023two} & \checkmark & & & P & E & & 4 & 16 & 50 \\ 
    \cite{mohri2024designing} & & \checkmark & & A & H & \checkmark & 1 & 26 & 296 \\ 
    \cite{lehmann2024matheuristic} & & & & A & H & & 6 & 4 & 100 \\ 
    \cite{achamrah2024gradient} & & \checkmark & & A & E & & 20 & 10 & 100 \\ 
    \cite{al2024multi} & \checkmark & & & A & H & \checkmark & 1 & 12 & 76 \\ 
    Our work & \checkmark & \checkmark & \checkmark & H & H & \checkmark & 10 & 109 & 1122 \\
    \bottomrule
  \end{tabular}
  \begin{tablenotes}
    \footnotesize
    \item[a] MP: Multi-period; CD: Containerized delivery; VR: Vehicle reposition.
    \item[b] Form.: Model formulation (A: Arc-based; P: Path-based; H: Hybrid arc-path-based); Sol.: Solution approach (H: Heuristic, E: Exact).
    \item[c] Real: Real dataset or case study.
  \end{tablenotes}
\end{threeparttable}
}
\end{table}

The two-echelon structure and the additional features introduced in the 2EVRP significantly increase the problem complexity. Given that the VRP is an NP-hard problem, the 2EVRP, as its generalization, is also NP-hard \citep{hemmelmayr2012adaptive}. Researchers have proposed various modeling formulations and solution approaches to address the 2EVRP. In terms of model formulation, three major paradigms have emerged. Arc-based formulations explicitly model decision variables for each arc in the logistic network, such as those by \cite{nguyen2012solving, yu2020two, lv2022two}. Path-based formulations, on the other hand, use variables corresponding to route selection, as seen in \cite{zeng2014hybrid, breunig2019electric, mohamed2023two}. More recently, a hybrid formulation integrating arc-based and path-based elements has been introduced in the 2EVRP by \cite{dellaert2021multi}. In their approach, the first-echelon delivery is modeled using an arc-based formulation, and the second-echelon routing is represented through a path-based modeling, with linking constraints to synchronize tours across both levels. While in this study, we propose a new hybrid formulation that integrates arc-based modeling for container routing with path-based modeling for vehicle service planning across both echelons. This formulation offers more detailed operational insights for practical applications but requires efficient solution approaches, particularly when applied to large-scale instances.

Regarding solution approaches of the 2EVRP, existing methods are broadly classified into exact methods (E) and heuristic approaches (H). Among these, exact methods, including dynamic programming \citep{baldacci2013exact}, branch-and-price \citep{dellaert2021multi}, and logic-based Benders decomposition \citep{mohamed2023two}, provide optimal solutions but are generally limited to small-scale problems due to computational complexity. To address more realistic settings, heuristic and metaheuristic methods have become the main solution approach in the 2EVRP literature. A mainstream approach relies on the large neighborhood search to explore the solution space and iteratively refine feasible solutions, as shown in \cite{hemmelmayr2012adaptive, breunig2019electric, lehmann2024matheuristic, muhlbauer2021parallelised}. Another common technique builds on echelon decomposition, initially proposed by \cite{crainic2011multi}, which separates the problem into two subproblems corresponding to the first and second echelons, and further decomposes the second-echelon problem into smaller VRPs for each satellite. \cite{mohri2024designing} further extended this heuristic by incorporating costs from both echelons when assigning customers to satellites. \cite{yu2020two} solved the first-echelon routing problem with exact methods, while constructing second-echelon routes heuristically with backtracking to ensure consistency. Building on prior studies, this work develops an echelon-region decomposition heuristic that incorporates a tailored capacity-aware flow assignment and rolling-horizon framework, capturing the hierarchical structure of the 2EVRP while balancing solution quality with computational efficiency.

Despite the development of advanced solution approaches, addressing real-world large-scale instances continues to be a major challenge. As illustrated in the Case study column of Table \ref{tab:lr_problem}, empirical studies on the 2EVRP are often limited in scope and scale. Most rely on synthetic datasets, typically restricted to single-depot settings with no more than 25 satellites and 300 customers. While \cite{lv2022two} addressed a relatively large dataset comprising 519 customers and 25 satellites, yet this scale remains insufficient to capture modern urban networks that may encompass thousands of customers. In this study, we advance the literature by solving a real-world 2EVRP instance comprising 10 depots, 109 satellites, and 1,122 customers, representing a scale that far exceeds that of most existing research. The findings of our case study aim to support the practical deployment of two-echelon delivery vehicle routing strategies in hyperconnected urban logistic networks.

\section{Problem description and formulation}
\label{sec:pdf}

\subsection{Problem description}

\begin{figure}[h!]
    \centering
    \includegraphics[width=15.8cm]{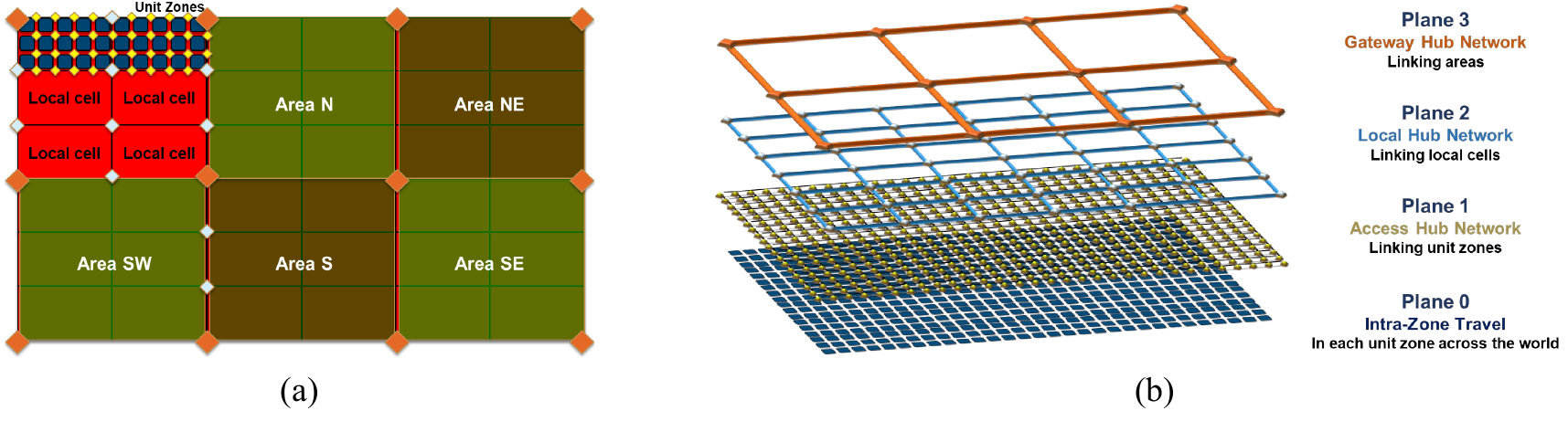} 
    \caption{{Structure of hyperconnected urban logistic networks \citep{montreuil2018urban}.}}
    \label{fig:hyperconnected_network}
\end{figure} 

This study introduces and addresses the dynamic two-echelon vehicle sharing and repositioning problem within a hyperconnected urban logistics system. The system is designed to enhance the operational efficiency and environmental sustainability of urban containerized deliveries through a PI-enabled three-tier hyperconnected logistic hub network. Under the PI paradigm, an urban territory is decomposed into unit zones, which are then aggregated into local cells, and these local cells are further grouped into larger urban areas, as illustrated in Figure \ref{fig:hyperconnected_network}. In a real-world context, a unit zone is the smallest modular service region, typically a residential block or industrial site; a local cell is a larger geographically coherent district; while an urban area usually represents a large segment of the city. Building upon this regional division, a three-tier interconnected meshed hub network is constructed, as depicted in Figure \ref{fig:hyperconnected_network}. The gateway hub network interconnects urban areas for managing shipments to and from external areas. Those gateway hubs are typically located at the edge of the city and are designed to handle large freight volumes and long-distance transfers, thereby facilitating the integration of the urban logistics system with the broader regional or national logistics network. Then the local hub network links adjacent local cells, with local hubs serving as satellites for intermediate freight transfer and coordination. Finally, the access hub network provides a finer level of granularity. Specifically, access hubs (e.g., automated parcel lockers), which are usually small and flexible, can be placed in a wide range of locations within communities for both customer pickup and courier intra-zone deliveries. 

\vspace{0.1in}
\begin{figure}[h!]
    \centering
    \includegraphics[width=16.8cm]{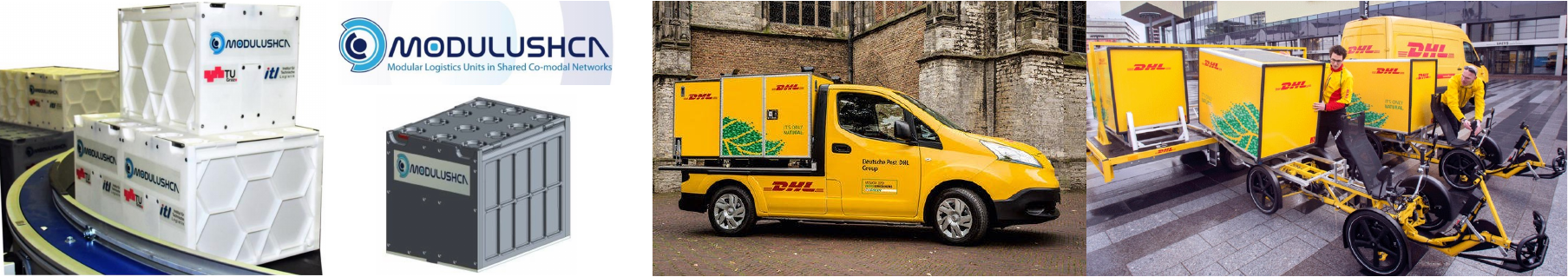} 
    \caption{{Illustrative examples of PI handling containers by \textit{Modulushca} project (www.modulshca.eu) and containerized last-mile delivery by \textit{DHL} (www.dhl.com).}}
    \label{fig:containerized_delivery}
\end{figure} 

Motivated by the concept of PI-containers \citep{montreuil2015modular}, we concentrate on containerized delivery operations across the three-tier hyperconnected network, assuming that all freight travels in standardized PI handling containers (H-containers). Compared with traditional pallets or cases, H-containers are lighter, stackable, and reusable, featuring a modular design that allows them to fit within delivery vehicles and supports automated loading and unloading systems. Their interlocking mechanism provides vertical stability during stacking, while maximizing spatial efficiency in storage and vehicle loading. A representative implementation of this concept is the \textit{Modulushca} project \citep{modulushca}, which developed a first-generation H-container prototype in practice, as shown in Figure \ref{fig:containerized_delivery}. By adopting such containers, vehicles swap full and empty containers without direct freight handling at each logistic hub, thereby reducing transfer time and enabling efficient operations. Once containers arrive at the access hubs, couriers unload packages to empty containers for subsequent swaps and may carry out final deliveries to customer-specified locations through intra-zone travel; however, this part of the operation is out of the scope of this work.

\begin{figure}[h!]
    \centering
    \includegraphics[width=10.5cm]{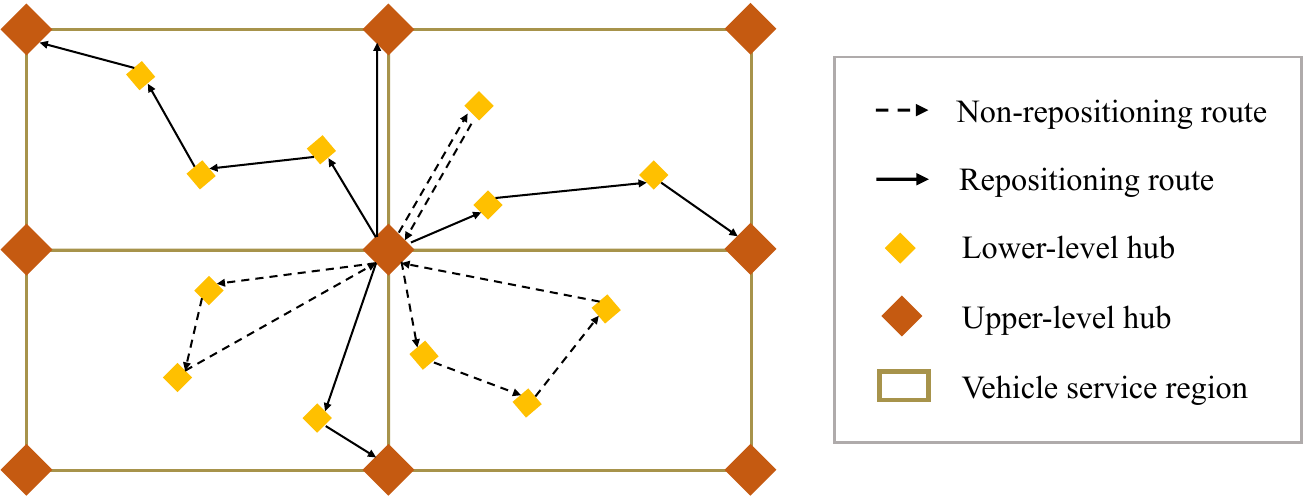}
    \caption{{Single-echelon vehicle route pattern diagram.}}
    \label{fig:vehicle_routes}
\end{figure}

In this study, we consider a two‑echelon vehicle routing structure for containerized deliveries, incorporating delivery vehicle sharing and repositioning across same‑level hubs. Both echelons operate based on predefined sets of potential scheduled service routes, as depicted conceptually in Figure \ref{fig:vehicle_routes}. In the first echelon, these routes connect gateway hubs with local hubs, while in the second echelon, they connect local hubs with access hubs. More specifically, delivery vehicles depart from upper-level hubs (i.e., gateway or local hubs), follow scheduled service routes to traverse one or multiple lower-level hubs (i.e., local or access hubs, respectively), and then either return to their origin hubs or end at different hubs for repositioning purposes. In addition, we further include the potential routes that traverse between upper-level hubs (e.g., gateway hubs), which enhances the connectivity among areas. Consequently, the number of vehicles at each gateway and local hub changes dynamically over time due to the repositions. We aim to determine which service routes to activate and how many delivery vehicles to assign to each active route in order to serve customer demand during each operational period. 

\vspace{0.1in}

Furthermore, when delivery demand surpasses the available capacity (i.e., the maximum number of available vehicles at hubs), outsourcing is allowed but penalized with an additional cost. This penalty can be interpreted as the direct cost of on-demand external fulfillment or the penalty of unmet customer demand. By adjusting this outsourcing penalty, logistics providers can balance investments in additional service frequencies with the flexibility of outsourced operations when demand surges, thereby accommodating different risk preferences and budget constraints. Moreover, the urban logistics environment is inherently dynamic, with demand for deliveries fluctuating significantly across different time slots of the day and days of the week. To capture these temporal dynamics and enable proactive service route planning, our problem is formulated in a multi-period setting. The planning horizon, such as several days or a week, is divided into several operational periods. This setting introduces additional complexity, as intertemporal dependencies arise between decision periods. Decisions made in one period, such as vehicle repositioning, directly impact the vehicle availability in subsequent periods.

\subsection{Formulation}
\label{sec:formulation}

To capture both service route planning and container routing, we formulate a multi-period arc- and path-based hybrid optimization model. At each period, the model determines: (i) the allocation of first‑ and second‑echelon vehicles to each potential service route (path-based decisions); (ii) the available number (i.e., inventory) of the first‑ and second‑echelon vehicles at gateway and local hubs; (iii) the routing of containers from origins to destinations using first‑ and second‑echelon deliveries (arc-based decisions); and (iv) the number of containers to be outsourced when service capacity is insufficient (arc-based decisions). The goal is to minimize total delivery operating costs while meeting the logistics needs of urban areas. The notation for the problem formulation is described as follows:

\begin{spacing}{1.28}
\begin{smallsize} 
\begin{longtable}[H]{ll}
\caption{{Notation for the multi-period arc- and path-based hybrid optimization model.}}\\
\midrule
\multicolumn{2}{l}{\textbf{Sets}} \vspace{-0.02in}\\
\midrule
$\boldsymbol{G}$ & Set of gateway hubs, indexed by $g$ \\
$\boldsymbol{L}$ & Set of local hubs, indexed by $l$ \\
$\boldsymbol{A}$ & Set of access hubs, indexed by $a$ \\
$\boldsymbol{R}$ & Set of FE routes, indexed by $r$ \\
$\boldsymbol{M}$ & Set of SE routes, indexed by $m$ \\
$\boldsymbol{T}$ & Set of time periods, indexed by $t$ \\
$\boldsymbol{A}^1$ & Set of arcs in the FE network, where $\boldsymbol{A}^1 = \{(i,j): i, j\in\boldsymbol{G}\cup\boldsymbol{L}\}$ \\
$\boldsymbol{A}^2$ & Set of arcs in the SE network, where $\boldsymbol{A}^2 = \{(i,j): i, j\in\boldsymbol{L}\cup\boldsymbol{A}\}$ \\
\midrule
\multicolumn{2}{l}{\textbf{Parameters}} \vspace{-0.02in} \\
\midrule
$A^S_{g r}$ (or $A^S_{l m}$) & Binary parameter, equals 1 if FE route $r$ (or SE route $m$) starts at hub $g$ (or $l$) \\
$A^E_{g r}$ (or $A^E_{l m}$) & Binary parameter, equals 1 if FE route $r$ (or SE route $m$) ends at hub $g$ (or $l$) \\
$B^{r}_{ij}$ & Binary parameter, equals 1 if arc $(i,j)\in \boldsymbol{A}^1$ is traversed by FE route $r$ \\
$B^{m}_{ij}$ & Binary parameter, equals 1 if arc $(i,j)\in \boldsymbol{A}^2$ is traversed by SE route $m$ \\
$D_{g a t}$ & Demand from gateway $g$ to access hub $a$ in period $t$ \\
$Q_1$ & Capacity of a FE vehicle \\
$Q_2$ & Capacity of a SE vehicle \\
$C^T_r$ & Scheduled service cost of using FE route $r$ \\
$C^T_m$ & Scheduled service cost of using SE route $m$ \\
$C^P_{i j}$ & On-demand penalty for transporting one container on arc $(i,j) \in \boldsymbol{A}^1 \cup \boldsymbol{A}^2$ \\
$I_g$ & Initial number of FE vehicles located at gateway hub $g$ \\
$I_l$ & Initial number of SE vehicles located at local hub $l$ \\
\midrule
\multicolumn{2}{l}{\textbf{Variables}} \vspace{-0.02in} \\
\midrule
$x_{r t}$ & Integer, number of FE vehicles assigned to route $r$ in period $t$ \\
$y_{m t}$ & Integer, number of SE vehicles assigned to route $m$ in period $t$ \\
$w_{g t}$ & Integer, number of available FE vehicles at gateway hub $g$ in period $t$ \\
$z_{l t}$ & Integer, number of available SE vehicles at local hub $l$ in period $t$ \\
$f_{ga,ij,t}$ & Integer, number of containers for demand from $g$ to $a$ transported on FE arc $(i,j)\in \boldsymbol{A}^1$ in period $t$ \\
$h_{ga,ij,t}$ & Integer, number of containers for demand from $g$ to $a$ transported on SE arc $(i,j)\in \boldsymbol{A}^2$ in period $t$ \\
$u_{i j t}$ & Integer, number of containers using on-demand service on FE arc $(i,j) \in \boldsymbol{A^1}$ in period $t$ \\
$v_{i j t}$ & Integer, number of containers using on-demand service on SE arc $(i,j) \in \boldsymbol{A^2}$ in period $t$ \\
\midrule
\end{longtable}
\end{smallsize}
\end{spacing}

Then the optimization model can be stated as follows:
\begin{midsize}
\begin{flalign}
& \min \quad \sum_{t \in \bm{T}} \Biggl[ 
\sum_{r \in \bm{R}} C^T_r \, x_{r t} + \sum_{m \in \bm{M}} C^T_m \, y_{m t}
+ \sum_{(i,j) \in \bm{A^1}} C^P_{i j} \, u_{i j t} + \sum_{(i,j) \in \bm{A^2}} C^P_{i j} \, v_{i j t} \Biggr]\label{eq:1} &
\end{flalign}
\vspace*{-0.24in}
\begin{flalign}
&s.t. \sum_{j:(i,j)\in \bm{A}^1} f_{ga,ij,t} + \sum_{j:(i,j)\in \bm{A}^2} h_{ga,ij,t} - \sum_{j:(j,i)\in \bm{A}^1} f_{ga,ji,t} - \sum_{j:(j,i)\in \bm{A}^2} h_{ga,ji,t} 
=\left\{\begin{array}{cl}
D_{gat}, & \text { if } i=g \\
-D_{gat}, & \text { if } i=a \\
0, & \text { otherwise}
\end{array} \right., &\notag \\
& \forall g\in \bm{G}, a\in \bm{A}, t\in \bm{T}, i \in \bm{G}\cup\bm{L}\cup\bm{A} \label{eq:2}
\end{flalign}
\vspace*{-0.4in}
\begin{flalign}
&\sum_{g\in \bm{G}} \sum_{a\in \bm{A}} f_{ga,ij,t} \leq \sum_{r\in \bm{R}} B^{r}_{ij}\, Q_1 \, x_{r t} + u_{i j t}, 
&\forall (i,j)\in \bm{A}^1, t\in \bm{T} \label{eq:3}\\
&\sum_{l\in \bm{L}} \sum_{a\in \bm{A}} h_{ga,ij,t} \leq \sum_{m\in \bm{M}} B^{m}_{ij}\, Q_2 \, y_{m t} + v_{i j t},
&\forall (i,j)\in \bm{A}^2, t\in \bm{T} \label{eq:4}\\
&\sum_{r\in \bm{R}} A^S_{g r}\, x_{r t} \leq w_{g t},
&\forall g\in \bm{G}, t\in \bm{T} \label{eq:5}\\
&\sum_{m\in \bm{M}} A^S_{l m}\, y_{m t} \leq z_{l t},
&\forall l\in \bm{L}, t\in \bm{T} \label{eq:6}\\
&w_{g, 0} = I_{g},
&\forall g\in \bm{G}\label{eq:7}\\
&w_{g t} = w_{g,t-1} - \sum_{r\in \bm{R}} A^S_{g r}\, x_{r,t-1} + \sum_{r\in \bm{R}} A^E_{g r}\, x_{r,t-1},
&\forall g\in \bm{G}, t\in \bm{T}\setminus\{0\}\label{eq:8}\\
&z_{l, 0} = I_{l},
&\forall l\in \bm{L} \label{eq:9}\\
&z_{l t} = z_{l,t-1} - \sum_{m\in \bm{M}} A^S_{l m}\, y_{m,t-1} + \sum_{m\in \bm{M}} A^E_{l m}\, y_{m,t-1},
&\forall l\in \bm{L}, t\in \bm{T}\setminus\{0\} \label{eq:10}\\
&x_{r t},\, y_{m t},\, w_{g t},\, z_{l t},\, f_{ga,ij,t},\, h_{ga,ij,t},\, u_{i j t},\, v_{i j t} \in \mathbb{Z}_{\geq 0}. \label{eq:11}
\end{flalign}
\end{midsize}

The objective function (\ref{eq:1}) minimizes the total operational costs, which include scheduled service operating costs and on-demand outsourcing penalties. Constraint (\ref{eq:2}) ensures unified container flow conservation across the entire two-echelon network for each origin-destination pair. Constraints (\ref{eq:3})-(\ref{eq:4}) enforce arc-level capacity limits in both echelons and ensure that, excluding the outsourcing part,
assigned vehicles are sufficient to satisfy flow requirements. Constraints (\ref{eq:5})-(\ref{eq:6}) govern node-level capacity limits at both gateway and local hubs. Constraints (\ref{eq:7})-(\ref{eq:8}) ensure the availability and assignment of first-echelon vehicles at gateway hubs over time, while constraints (\ref{eq:9})-(\ref{eq:10}) impose similar inventory constraints for second-echelon vehicles at local hubs. Lastly, constraints (\ref{eq:11}) specify the domains of the decision variables.

\section{Solution approach}
\label{sec:sa}

Our proposed heuristic algorithm addresses the large-scale 2EVRP by integrating (i) capacity-aware flow assignment, (ii) echelon-region decomposition-based route optimization, and (iii) a rolling-horizon framework. Specifically, we first estimate per-unit route costs for both echelons (Section \ref{sec:route-cost-estimation}) and then allocate demand flow to local hubs via a lightweight capacity-aware flow optimization (Section \ref{sec:capacity-aware-opt}). Conditioning on these flows, the first-echelon and second-echelon routing problems can be decoupled: we solve a single first-echelon subproblem over the network of gateway and local hubs, and a set of second-echelon subproblems, one per urban area, each solved independently on its sub-network (Section \ref{sec:decomp}). Finally, we solve the problem within a rolling-horizon framework to incorporate forward-looking in vehicle repositioning strategies. Moreover, this enables frequent updates of hub capacities (i.e., vehicle availability), leading to better capacity-aware flow assignment solutions (Section \ref{sec:rolling}).

\subsection{Capacity-aware flow assignment}
\label{sec:capacity-aware-assignment}
\subsubsection{Route cost estimation}
\label{sec:route-cost-estimation}

We require a fast and robust proxy for the marginal first-echelon price of sending one unit from a gateway hub $g\in\bm{G}$ to a local hub $l\in\bm{L}$. Let $\bm{R}_g$ denote the set of candidate first-echelon routes starting from gateway hub $g$. Each route $r\in\bm{R}_g$ bears a per-unit service cost $C_r>0$, computed as the sum of an operating cost related to travel time and a repositioning surcharge, i.e.,
\begin{equation}
C_r \;=\; C_{\mathrm{op/hr}} \cdot \text{TravelTime}_r  \;+\; C_{\mathrm{reposition}}\cdot \mathbf{1}\{\text{reposition}_r\}.
\end{equation}

For a fixed pair $(g,l)\in\bm{G}\times\bm{L}$, we define $\lambda_{g l}$ as the optimal value of a unit-flow linear program (LP) that selects the cheapest route in $\bm{R}_g$ that visits local hub $l$:

\begin{equation}
\label{eq:unitflow}
\begin{aligned}
\lambda_{g l} \;=\; \min \quad & \sum_{r\in\bm{R}_g} C_r\, x_r \\
\text{s.t.}\quad & \sum_{r\in\bm{R}_g} B^r_l\, x_r \;\ge\; 1,\\
\quad & x_r \;\ge\; 0,
\end{aligned}
\end{equation}
where $B^r_l$ is a binary parameter to indicate the coverage of local hubs by routes. $B^r_l=1$ if and only if local hub $l$ is visited by route $r$, and $B^r_l=0$ otherwise. In implementation, we construct one unit-flow LP per gateway hub $g$ containing the variables $\{x_r\}_{r\in\bm{R}_g}$, a fixed linear objective $\sum_r C_r x_r$, and a cover constraint whose coefficients $\{B^r_l\}_{r\in\bm{R}_g}$ are edited in place for each target local hub $l$. Solving~\eqref{eq:unitflow} sequentially over all $l\in\bm{L}$ while reusing the same model preserves the simplex basis across re-optimizations. Therefore, employing dual simplex delivers warm-started solves that converge in a few iterations per local hub. If no route in $\bm{R}_g$ covers $l$, we set $\lambda_{g l}=+\infty$.

Since the first-echelon delivery allows freight to first transit from one gateway hub to another gateway hub (e.g., $g\to g'$) and then transfer down to the local hub layer, we recompute $\lambda_{g l}$ by including inter-gateway hub connections. Let $\mu_{g g'}$ denote the unit price to reach gateway hub $g'$ from gateway hub $g$, obtained by solving the same unit-flow model~\eqref{eq:unitflow} with a coverage binary parameter $B^r_{g'}=1$ if route $r$ terminates at gateway hub $g'$. We then compute the all-pairs transfer $\mu^\star$ by updating $\mu$ via the Floyd-Warshall algorithm:
\begin{equation}
\mu^\star_{i j} \;\leftarrow\; \min\; \bigl\{\mu_{i j},\;\mu^\star_{i k}+\mu^\star_{k j}\bigr\}, \qquad \forall\, i,j,k\in\bm{G}
\end{equation}
which captures the cheapest path between any pair of gateway hubs. Consequently, the per-unit marginal price $\lambda_{g l}$ can be updated by allowing intermediate gateway hubs as follows:
\begin{equation}
\label{eq:closure}
\lambda_{g l} \;\leftarrow\; \min\; \Bigl\{\,\lambda_{g l},\;\min_{g'\in\bm{G}}\bigl(\mu^\star_{g g'}+\lambda_{g' l}\bigr)\Bigr\},
\end{equation}
so that if reaching $l$ is cheaper via $g'$ than directly from $g$ (e.g., $\lambda_{g l}=+\infty$), the updated value reflects this multi-route path. Following a similar approach, we further compute the marginal second-echelon price of sending one unit from a local hub $l\in\bm{L}$ to an access hub $a\in\bm{A}$, referred to as $\lambda_{l a}$. 

\subsubsection{Capacity-aware flow optimization}
\label{sec:capacity-aware-opt}

Given the marginal prices $\lambda_{g l}$ and $\lambda_{l a}$ obtained in the previous subsection, we solve a capacitated flow assignment LP to determine how to route the splittable container flow $D_{g a t}$ through local hub $l$. Let $s_{g l a t} \in \mathbb{Z}_{\geq 0}$ denotes the assigned container flow from gateway hub $g$ to access hub $a$ in period $t$ transiting through local hub $l$. The target is to minimize the total costs across both echelons subject to demand conservation and hub capacity (i.e., the number of available vehicles at hubs). However, in practice, infeasibility may arise due to temporary capacity shortages. To address this, we further modify the capacity restriction as soft constraints by introducing slack variables $e^1_{g t}$ and $e^2_{l t}$, where the variables are associated with positive penalties $C^P_g$ and $C^P_l$, respectively. Then, the capacitated flow optimization problem can be formulated as:

\begin{equation}
\label{eq:assign}
\begin{aligned}
z^*(t, \hat{w}_{gt}, \hat{z}_{lt}) \;=\; \min \quad & \sum_{g\in\bm{G}}\sum_{a\in\bm{A}}\sum_{l\in\bm{L}}
s_{g l a t}\bigl(\lambda_{g l}+\lambda_{l a}\bigr) +\sum_{g\in\bm{G}}C^P_g e^1_{g t} + \sum_{l\in\bm{L}}C^P_l e^2_{l t} \\
\text{s.t.}\quad &
\sum_{l\in\bm{L}} s_{g l a t} \;=\; D_{g a t},\quad \forall g\in\bm{G}, a\in\bm{A}\\
\quad & \sum_{l\in\bm{L}}\sum_{a\in\bm{A}} s_{g l a t} \;\le\; Q_1 \hat{w}_{g t} + e^1_{g t}, \quad \forall g\in\bm{G}\\
\quad & \sum_{g\in\bm{G}}\sum_{a\in\bm{A}} s_{g l a t} \;\le\; Q_2 \hat{z}_{l t} + e^2_{l t}, \quad \forall l\in\bm{L}\\
& s_{g l a t}, e^1_{g t}, e^2_{l t} \in \mathbb{Z}_{\geq 0},
\end{aligned}
\vspace{0.1in}
\end{equation}
where $w_{g t}$ and $z_{l t}$ denote the available number of vehicles at hub $g$ and $l$ in period $t$ as we introduced earlier.  

For each period $t$, we solve problem~\eqref{eq:assign} given $\hat{w}_{gt}$ and $\hat{z}_{lt}$, derived from the solutions in the previous period. Then, the final assignment yields the aggregated flows:

\vspace{-0.05in}
\begin{equation}
\label{eq:aggflows}
D_{g l t}\;=\;\sum_{a\in\bm{A}} s_{g l a t},
\qquad
D_{l a t}\;=\;\sum_{g\in\bm{G}} s_{g l a t},
\end{equation}
which provide inputs to the first-echelon and second-echelon subproblems. Overall, our approach accommodates splittable deliveries and captures vehicle availability without introducing too much complexity in the flow assignment step.

\subsection{Heuristic echelon-region decomposition}
\label{sec:decomp}
Given the flow assignment $D_{g l t}$ and $D_{l a t}$ from Section \ref{sec:capacity-aware-assignment}, we decompose our 2EVRP problem (\ref{eq:1})-(\ref{eq:11}) into (i) a first-echelon optimization over the gateway-local hub network and (ii) a set of independent second-echelon optimizations, each restricted to a single urban area, over the local-access hub network. We follow the notation that is introduced in Section \ref{sec:formulation}. For period $t$, given $D_{g l t}$ specifies the number of containers to be shipped from gateway hub $g$ to local hub $l$ and $\hat{w}_{gt}$ denotes the number of available delivery trucks at hub $g$, the first-echelon model at period $t$ can be expressed as follows: 

\vspace{-0.2in}
\begin{midsize}
\begin{flalign}
z^{1*}(t, D_{g l t}, \hat{w}_{g t}) \;=\;  & \min \quad 
\sum_{r \in \bm{R}} C^T_r \, x_{r t} + \sum_{(i,j) \in \bm{A^1}} C^P_{i j} \, u_{i j t}\\
&s.t. \sum_{j:(i,j)\in \bm{A}^1} f_{gl,ij,t} - \sum_{j:(j,i)\in \bm{A}^1} f_{gl,ji,t} 
=\left\{\begin{array}{cl}
D_{glt}, & \text { if } i=g \\
-D_{glt}, & \text { if } i=l \\
0, & \text { otherwise}
\end{array} \right.,  &\hspace{-0.1in}\forall g\in \bm{G}, l\in \bm{L}, i \in \bm{G}\cup\bm{L} \label{eq:111}\\
&\sum_{g\in \bm{G}} \sum_{a\in \bm{A}} f_{ga,ij,t} \leq \sum_{r\in \bm{R}} B^{r}_{ij}\, Q_1 \, x_{r t} + u_{i j t}, 
\quad &\forall (i,j)\in \bm{A}^1 \label{eq:112}\\
&\sum_{r\in \bm{R}} A^S_{g r}\, x_{r t} \leq \hat{w}_{g t}, \quad  &\forall g\in \bm{G} \label{eq:113}\\
&x_{r t},\, f_{ga,ij,t},\, u_{i j t} \in \mathbb{Z}_{\geq 0}.
\end{flalign}
\end{midsize}

Correspondingly, for each period $t$, given demand flow $D_{lat}$ and local hub capacity $\hat{z}_{lt}$, the second-echelon optimization model can be stated as:
\begin{midsize}
\begin{flalign}
z^{2*}(t,D_{l a t}, \hat{z}_{l \tau}) \;=\; & \min \quad 
\sum_{m \in \bm{M}} C^T_m \, y_{m t}
+ \sum_{(i,j) \in \bm{A^2}} C^P_{i j} \, v_{i j t} \label{eq:210}\\
& s.t. \sum_{j:(i,j)\in \bm{A}^2} h_{la,ij,t} - \sum_{j:(j,i)\in \bm{A}^2} h_{la,ji,t} 
=\left\{\begin{array}{cl}
D_{lat}, & \text { if } i=l \\
-D_{lat}, & \text { if } i=a \\
0, & \text { otherwise}
\end{array} \right.,
&\hspace{-0.1in}\forall l\in \bm{L}, a\in \bm{A}, i \in \bm{L}\cup\bm{A} \label{eq:211}\\
&\sum_{l\in \bm{L}} \sum_{a\in \bm{A}} h_{ga,ij,t} \leq \sum_{m\in \bm{M}} B^{m}_{ij}\, Q_2 \, y_{m t} + v_{i j t},
\quad &\forall (i,j)\in \bm{A}^2 \label{eq:212}\\
&\sum_{m\in \bm{M}} A^S_{l m}\, y_{m t} \leq \hat{z}_{l t},\quad &\forall l\in \bm{L}, t \in \bm{T}_\tau\\
&y_{m t},\, h_{ga,ij,t},\, v_{i j t} \in \mathbb{Z}_{\geq 0} \label{eq:213}.
\end{flalign}
\end{midsize}

Due to the substantial number of local and access hubs, the second-echelon subproblem continues to exhibit high computational complexity. To further improve the algorithm efficiency, we exploit the spatial structure of the second echelon delivery in the hyperconnected network. Specifically, let $\bm{\Gamma}$ represent the set of urban areas and, for each $\gamma\in\bm{\Gamma}$, $\bm{L}_\gamma\subseteq\bm{L}$ and $\bm{A}_\gamma\subseteq\bm{A}$ denote the local and access hubs in urban area $\gamma$. We induce the sub-network $\bm{A}^2_\gamma$ and the route set $\bm{M}_\gamma$ by restricting arcs and routes to those whose start and end points lie within $\bm{L}_\gamma\cup\bm{A}_\gamma$. Similarly, the demand $D_{lat}$ is partitioned by $\gamma$, decomposing model~\eqref{eq:210}-\eqref{eq:213} into $|\bm{\Gamma}|$ independent MIPs that can be solved in parallel. In the implementation, this decomposition approach significantly reduces model size and solution time.

\subsection{Forward-looking rolling horizon framework}
\label{sec:rolling}
We solve the decomposed models within a rolling-horizon framework to incorporate forward-looking in decision-making while keeping the relatively frequent update of vehicle availability for capacity-aware flow optimization. In practice, if the sliding window is too short (e.g., solving only one period at a time), the model captures limited temporal interactions and offers fewer opportunities for optimizing vehicle repositioning across periods. Conversely, if the window length is long (e.g., spanning the full planning horizon), hub capacity updates become infrequent, reducing the accuracy of the flow assignment step and potentially degrading solution quality. To balance this trade-off, we adopt an intermediate window size that allows for sufficient look-ahead while maintaining frequent hub capacity updates over the planning horizon.

\vspace{0.1in}
\begin{figure}[h!]
    \centering
    \includegraphics[width=13.6cm]{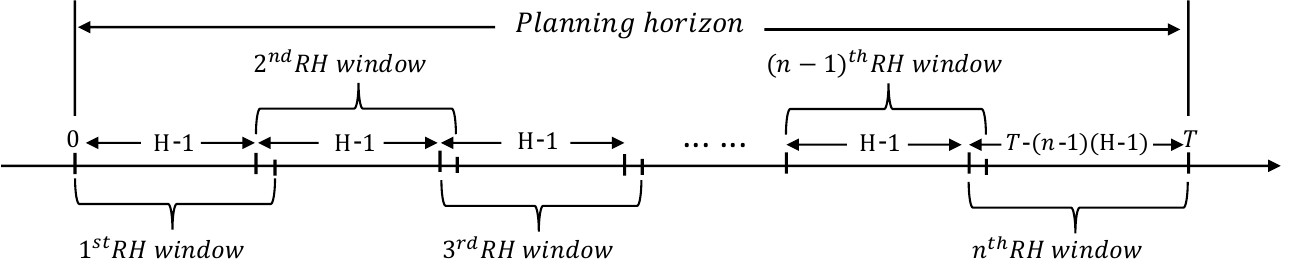}
    \caption{{An illustration of rolling horizon scheduling.}}
    \label{fig:rolling_horizon}
\end{figure} 

Figure \ref{fig:rolling_horizon} illustrates the forward-looking rolling horizon approach. Specifically, the planning horizon $\bm{T} = \{0, 1, ...,T\}$ can be discretized into $|\bm{T}|$ periods. Let $H\in\mathbb{Z}_{>2}$ denote the window length. For a window starting at period $\tau$, we form $\bm{T}_\tau=\{\tau,\ldots,\tau+H-1\}\cap\bm{T}$ and define the implementation set $\bm{I}_\tau\subseteq\bm{T}_\tau$ as all but the last period in the window, i.e., $\bm{I}_\tau=\{\tau,\ldots,\tau+H-2\}$. Only when the window reaches the end of the planning horizon, we have $\bm{I}_\tau=\bm{T}_\tau$. Moreover, we initialize the pre-decision vehicle availability by $\{w_{g\tau}\}_{g\in\bm{G}}$ and $\{z_{l\tau}\}_{l\in\bm{L}}$ computed from the previous implementation. Given these states, we then proceed in three steps: First, we solve the capacity-aware flow assignment~\eqref{eq:assign} for each period $t\in\bm{T}_\tau$ using the per-unit prices $(\lambda_{g l}+\lambda_{l a})$ and bounding the total assigned flow based on the vehicle numbers $w_{g\tau}$ and $z_{l\tau}$. This yields period-wise aggregate flows $\{D_{glt}\}_{t\in\bm{T}_\tau}$ and $\{D_{lat}\}_{t\in\bm{T}_\tau}$. Second, conditioning on these flows, we optimize the first-echelon model over all $t\in\bm{T}_\tau$. Third, we solve the second-echelon models independently for each urban area $\gamma \in \bm{\Gamma}$, on its induced sub-network over the same $\bm{T}_\tau$. 

We remark that, in the first step, we solve the single-period flow assignment problem \eqref{eq:assign} using the same capacity limits $w_{g\tau}$ and $z_{l\tau}$ for each period $t\in\bm{T}_\tau$. This step estimates transfer flows at local hubs under the assumption that vehicle repositioning is not enabled. In the second step, to capture potential savings from repositioning opportunities, we re-optimize the decomposed models across multiple periods to account for these reposition benefits. While this approach does not guarantee optimality, our empirical results show that it consistently yields high-quality solutions while ensuring solvability for large-scale instances. It is worth noting that decomposed models are presented as single-period formulations in Section \ref{sec:decomp}. To enable cross-period optimization within the rolling-horizon framework, we extend them to solve over a set of periods as follows:

\vspace{-0.15in}
\begin{midsize}
\begin{equation}
\label{eq:fe}
\begin{aligned}
z^{1*}(\bm{T}_\tau, D_{g l t}, \hat{w}_{g \tau}) \;=\;  & \min \quad 
\sum_{t \in \bm{T}_\tau} \left[ \sum_{r \in \bm{R}} C^T_r \, x_{r t} + \sum_{(i,j) \in \bm{A^1}} C^P_{i j} \, u_{i j t} \right]\\
&s.t. (\ref{eq:111}), (\ref{eq:112}) \text{ for } \forall t \in \bm{T}_\tau\\
&\sum_{r\in \bm{R}} A^S_{g r}\, x_{r t} \leq w_{g t}, \quad  \forall g\in \bm{G}, t \in \bm{T}_\tau\\
&w_{g, \tau} = \hat{w}_{g \tau}, \quad \forall g\in \bm{G}\\
&w_{g t} = w_{g,t-1} - \sum_{r\in \bm{R}} A^S_{g r}\, x_{r,t-1} + \sum_{r\in \bm{R}} A^E_{g r}\, x_{r,t-1}, \quad \forall g\in \bm{G}, t\in \bm{T}_\tau\setminus\{\tau\}\\
&x_{r t},\, f_{ga,ij,t},\, w_{g t}, u_{i j t} \in \mathbb{Z}_{\geq 0}.
\end{aligned}
\end{equation}

\quad
\begin{equation}
\label{eq:se}
\begin{aligned}
z^{2*}(\bm{T}_\tau,D_{l a t}, \hat{z}_{l \tau}) \;=\; & \min \quad 
\sum_{t \in \bm{T}_\tau} \left[ \sum_{m \in \bm{M}} C^T_m \, y_{m t}
+ \sum_{(i,j) \in \bm{A^2}} C^P_{i j} \, v_{i j t} \right]\\
& s.t. (\ref{eq:211}), (\ref{eq:212}) \text{ for } \forall t \in \bm{T}_\tau\\
&\sum_{m\in \bm{M}} A^S_{l m}\, y_{m t} \leq z_{l t},\quad \forall l\in \bm{L}, t \in \bm{T}_\tau\\
&z_{l, \tau} = \hat{z}_{l \tau}, \quad \forall l\in \bm{L} \\
&z_{l t} = z_{l,t-1} - \sum_{m\in \bm{M}} A^S_{l m}\, y_{m,t-1} + \sum_{m\in \bm{M}} A^E_{l m}\, y_{m,t-1}, \quad \forall l\in \bm{L}, t\in \bm{T}_\tau\setminus\{\tau\}\\
&y_{m t},\, h_{ga,ij,t},\, z_{lt}, v_{i j t} \in \mathbb{Z}_{\geq 0}.
\end{aligned}
\end{equation}
\end{midsize}

\vspace{0.05in}
This modification restores variables $w_{gt}$ and $z_{lt}$ along with constraints \eqref{eq:8} and \eqref{eq:10} to update vehicle availability across periods, with the first period initialized by $w_{g\tau}$ and $z_{l\tau}$, respectively. Finally, after solving the optimization models, we implement the decisions for all $t\in\bm{I}_\tau$: scheduled service decisions $(x^*, y^*)$ and on-demand service usages $(u^*, v^*)$ are fixed, container flow strategies $(f^*, h^*)$ are determined, vehicle availability $(w^*, z^*)$ is updated, and period costs are recorded. The planning horizon then advances by $H-1$ periods, and the procedure repeats until all periods have been executed. The overall heuristic framework is summarized in Algorithm \ref{alg:heuristic}.

\vspace{-0.05in}
\begin{algorithm}[h!]
\begin{midsize} 
\begin{spacing}{0.9}
\caption{{Echelon-region decomposition heuristic for 2EVRP}}
\vspace{0.05in}
\KwData{Initial vehicle number $I_g,I_l$, marginal prices $\{\lambda_{g l}\},\{\lambda_{l a}\}$, period set $\bm{T}$, window length $H$}
Sliding window start period $\tau \gets 0$\;
Vehicle numbers $w_{g0}\gets I_g,~z_{l0}\gets I_l$\;  Objective value $z^* \gets 0$\;
\vspace{0.05in}

\While{$\max\{\bm{T}_\tau\}<T$}{
  Update sliding window $\bm{T}_\tau \gets \{\tau,...,\tau+H-1\}\cap\bm{T}$ and implementation set $\bm{I}_\tau \gets \{\tau,...,\tau+H-2\}\cap\bm{T}$\;
  Update $\hat{w}_{g, \tau}$ and $\hat{z}_{l, \tau}$\;
  \vspace{0.1in}
  
  \Comment*[h]{/* Capacity-aware flow assignment */}\\
  Solve the problem \eqref{eq:assign} for each period $t \in \bm{T}_\tau$ with inputs $(\hat{w}_{g \tau},\hat{z}_{l \tau})$ and marginal prices $(\lambda_{g l},\lambda_{l a})$\; 
  Compute $D_{g l t}$ and $D_{l a t}$ for each period $t \in \bm{T}_\tau$ based on equation \eqref{eq:aggflows}\; 
  \vspace{0.1in}
  
  \Comment*[h]{/* Decompose problem by echelon and region */}\\
  Solve the problem \eqref{eq:fe} over $\bm{T}_\tau$ with input $D_{g l t}$ and $\hat{w}_{g \tau}$\; 
  Record the optimal solution $\{x^*_t,f^*_t,w^*_t,u^*_t\}_{t \in \bm{I}_\tau}$ and optimal value $z^{1*}(\bm{I}_\tau)$\;
  \For{urban area $\gamma\in\bm{\Gamma}$}{
  Restrict network set to that area and solve the problem \eqref{eq:se} over $\bm{T}_\tau$ with input $D_{l a t}$ and $\hat{z}_{l \tau}$\; 
  Record the per-area optimal solution $\{y^*_t(\gamma), v^*_t(\gamma), z^*_t(\gamma), h^*_t(\gamma)\}_{t \in \bm{I}_\tau}$ and optimal value $z^{2*}(\bm{I}_\tau, \gamma)$\; 
  }
  Calculate $z^{2*}(\bm{I}_\tau) \gets \sum_{\gamma\in\bm{\Gamma}} z^{2*}(\bm{I}_\tau,\gamma)$\;
  Update total objective value $z^* \gets z^*+(z^{1*}(\bm{I}_\tau)+z^{2*}(\bm{I}_\tau))$\;
  \vspace{0.1in}

  \Comment*[h]{/* Roll the horizon */}\\
  $\tau \gets \tau+H-2$\;
}
\label{alg:heuristic}
\end{spacing}
\end{midsize}
\vspace{-0.05in}
\end{algorithm}

\section{Case study}
\label{sec:cs}
\subsection{Dataset and parameter settings}
\label{sec:dps}
\begin{figure}[h!]
    \centering
    \includegraphics[width=16.1cm]{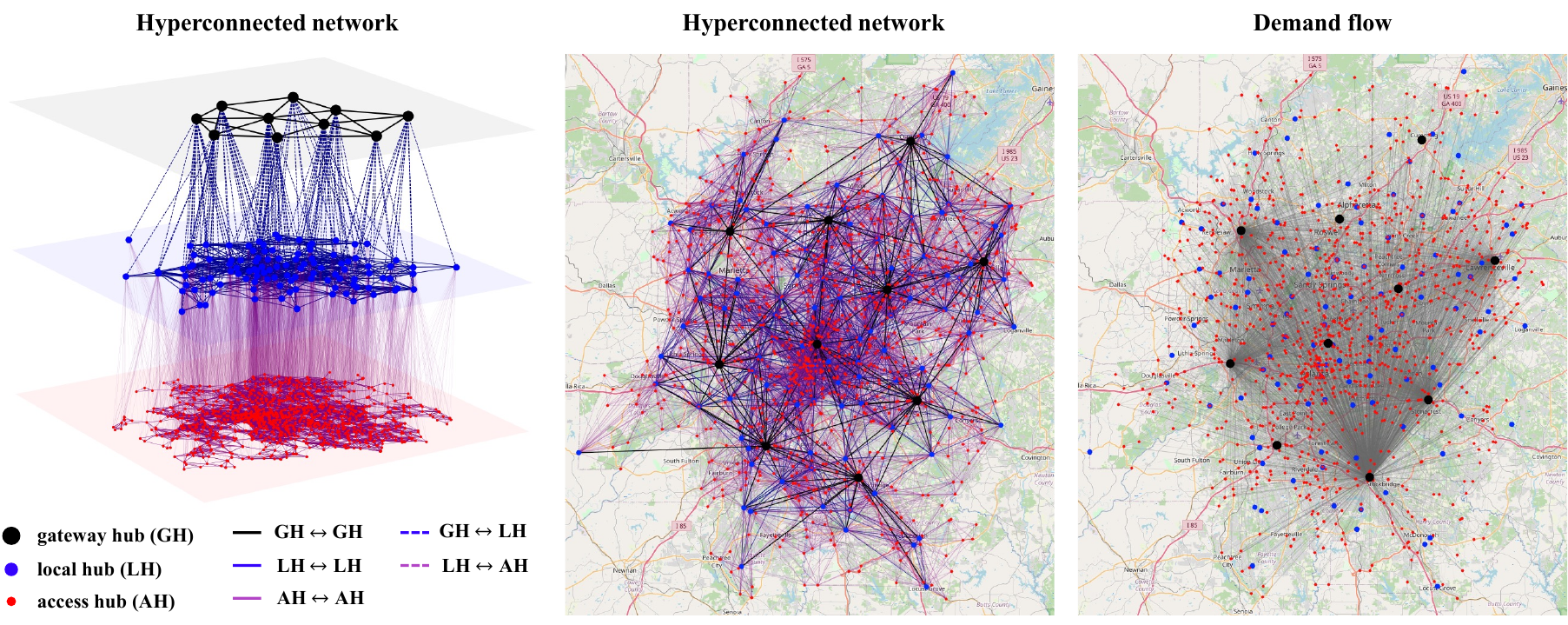}
    \caption{{Hub locations, service region mapping, and three-tier hyperconnected network \citep{liu2025network}.}}
    \label{fig:network_info}
\end{figure}

In this study, we consider a setting where logistics providers aim to deliver freight from large distribution centers (i.e., gateway hubs) to locker banks near customers (i.e., access hubs) using a containerized delivery in urban areas. We select the Atlanta metropolitan area as our research testbed, representing a large, complex urban environment with significant freight activity. Figure \ref{fig:network_info} presents the three-tier hyperconnected logistic hub network for this area. The hub locations are identified from the CoStar database \citep{costar}, a source of commercial real estate information, ensuring that the selected sites correspond to actual warehouses and distribution centers suitable for logistics operations. This extensive network, which forms the basis for our experiments, comprises a total of 10 gateway hubs, 109 local hubs, and 1,122 access hubs. In addition, as illustrated in the multi-layer visualization, this hyperconnected network features both vertical connectivity between tiers (e.g., gateway to local hubs, local to access hubs) and horizontal connectivity within tiers (e.g., between local hubs), enabling flexible freight consolidation and routing. Figure \ref{fig:network_info} further plots the origin-destination demand pairs within the network, visualizing the demand flow from gateway hubs (black points) to access hubs (red points) as grey lines, representing the primary delivery tasks. This demand flow is estimated based on the Freight Analysis Framework (FAF) database \citep{faf}, which provides comprehensive estimates of freight movement across the United States. It is important to note that in this study, we take the network design and OD flows as given and focus on two-echelon scheduled service route planning to serve customer demand; while the detailed information about network design and flow estimation was established in our prior research \citep{liu2025network}.

The longest planning horizon under consideration in our experiments spans 9 periods, equivalent to 3 days, with each day divided into three operational periods. This setting reflects a realistic operational context where a logistics provider forecasts demand up to three days ahead to plan its service capacity and resource allocations while retaining operational agility. In our experiment section, we test the proposed solution approach on varying planning horizons (3, 6, and 9 periods). Increasing the horizon grows the size of the problem instance. This enables us to study the scalability of our methods under increasing complexity and to observe how horizon depth influences vehicle repositioning among periods. Fleet sizing follows a two-echelon design: we assume that each first-echelon delivery truck has a capacity of 10 standardized containers, while each second-echelon delivery van can carry up to 4 containers. Moreover, each container is designed to hold a maximum of 1.5 tons of freight. These capacity parameters align with the medium- and light-duty vehicles used in urban logistics but can be easily adjusted for other city contexts or vehicle classes with minimal model changes. The scheduled service operating costs are estimated based on realistic factors, including fleet costs, driver wages, vehicle speeds, and travel distances. All of these values are calibrated to contemporary urban operations. The on-demand outsourcing cost is set at a high level to model it as a premium and less desirable option used only to handle exceptional demand surges that exceed the network's capacity. Lastly, the environmental parameters, including CO$_2$-eq emissions and energy consumption for each fuel technology, are derived from the R\&D GREET.Net model (Software v1.3.0.14345) \citep{GREET2025Updates}. A complete list of the corresponding parameter values used in the study is detailed in \ref{appendix_a}.

To generate our candidate route sets, we employ a two-stage approach that begins by enumerating a full set of potential paths based on predefined structural patterns for each echelon. This large initial set is then pruned using a pre-processing procedure, which filters out inefficient routes containing segments that exceed specified distance thresholds. Finally, a coverage restoration is applied to ensure network integrity by verifying that every access hub can be served from any gateway hub within the network, reintroducing the routes if any were removed during pruning. This process effectively creates a computationally manageable yet comprehensive set of routes for the optimization model. The details of this route generation are outlined in \ref{appendix_b}.

In this study, we design a comprehensive set of numerical experiments. In Section \ref{sec:cp_results}, a computational study is presented to assess the scalability of our proposed heuristic across various instance sizes. Section \ref{sec:ne} evaluates the effectiveness of the two-echelon (i.e., three-tier) hyperconnected delivery system by comparing its performance against a single-echelon (i.e., two-tier) counterpart under different market share scenarios. In Section \ref{sec:vsr}, we quantify the benefits derived from the vehicle sharing and repositioning strategies by comparing system performance with and without this operational flexibility. Subsequently, Section \ref{sec:ei} investigates the environmental impact of the urban delivery system, where we test vehicle fleets with different fuel types. Finally, in Section \ref{sec:sen_analysis}, we conduct a sensitivity analysis on key cost parameters to evaluate the robustness of the model's solutions and strategic decisions.

\subsection{Computational study}
\label{sec:cp_results}

To evaluate the performance and scalability of the proposed solution approach, we construct two region instances with varying network sizes and complexities. As detailed in Table \ref{tab:instance}, Region 1 represents a compact sub-network encompassing three urban areas with an average daily flow of 3,260 containers. While Region 2 scales this up to a large-scale network covering all ten urban areas shown in Figure \ref{fig:network_info}. This region features a significantly more complex structure with an average daily flow of 105,141 containers, thereby providing a challenging testbed for the tractability of the solution approaches.

\vspace{0.1in}
\begin{table}[h!]
\begin{center}
\caption{{Region instance characteristics.}}
\resizebox{0.73\linewidth}{!}{
\begin{threeparttable}
\begin{tabular}{cccccccccc}
\toprule
 & \#GHs & \#LHs & \#AHs & \#FE routes & \#SE routes & Daily flow in tons\\
\midrule
\textbf{\small{Region 1}} & 3 & 32 & 301 & 122 & 1,516 & 4,898 \\
\textbf{\small{Region 2}} & 10 & 109 & 1,122 & 503 & 10,331 & 157,230 \\ 
\bottomrule
\end{tabular}
\end{threeparttable}
}
\label{tab:instance}
\end{center}
\vspace{-6pt}
\end{table}

The computational experiment is designed to compare the performance of four distinct solution strategies across instances of varying sizes. The strategies are: (i) a Gurobi MILP solve (Gurobi), serving as a benchmark for solution quality; (ii) a baseline decomposition strategy (Decomp heuristic), where the flow assignment step assigns the LH closest to the destination as a satellite; (iii) the decomposition strategy with a rolling horizon framework (Decomp heuristic + RH), where the flow assignment step assigns the LH closest to the destination as a satellite; and (iv) our proposed enhanced decomposition strategy integrated with a rolling horizon framework and capacity-aware flow assignment (Enhanced decomp heuristic + RH). The instance sizes are varied along two dimensions: 1) the network size (i.e., Region 1 and Region and 2) the length of the planning horizon, which is set to 3, 6, and 9 periods, corresponding to one, two, and three days of operation, respectively. Moreover, we evaluate the method performance based on the following key metrics: the final objective function value in dollars (Objective), the final MIP optimality gap (\%Gap), the relative gap compared to the Gurobi benchmark solution (\%Rel Gap), and the total computational runtime in hours (Runtime).

\vspace{0.1in}
\begin{table}[h!]
\begin{center}
\caption{{Computational performance comparison across different solution approaches.}}
\resizebox{0.95\linewidth}{!}{
\begin{threeparttable}
\begin{tabular}{cccccccccc}
\toprule
 & & \multicolumn{4}{c}{\textbf{Gurobi}} & \multicolumn{4}{c}{\textbf{Decomp heuristic}} \\
\cmidrule(lr){3-6} \cmidrule(lr){7-10}
Region & \#Periods & Objective & \%Gap & \%Rel Gap & Runtime & Objective & \%Gap & \%Rel Gap & Runtime \\
\midrule
\textbf{Region 1} & 3 & 17,094 & 1.18 & 0.00 & 12.00 & 23,512 & 0.83 & 27.30 & 0.12 \\
         & 6 & 31,866 & 1.83 & 0.00 & 12.00 & 43,728 & 0.82 & 27.13 & 0.26 \\
         & 9 & 44,813 & 1.90 & 0.00 & 12.00 & 61,743 & 0.82 & 27.42 & 0.41 \\
\midrule
 & & \multicolumn{4}{c}{\textbf{Decomp heuristic + RH}} & \multicolumn{4}{c}{\textbf{Enhanced decomp heuristic + RH}} \\
\cmidrule(lr){3-6} \cmidrule(lr){7-10}
Region & \#Periods & Objective & \%Gap & \%Rel Gap & Runtime & Objective & \%Gap & \%Rel Gap & Runtime \\
\midrule
\textbf{Region 1} & 3 & 23,506 & 0.90 & 27.28 & 0.12 & 17,477 & 0.70 & 2.19 & 0.11 \\
         & 6 & 43,732 & 0.88 & 27.13 & 0.33 & 32,512 & 0.70 & 1.99 & 0.30 \\
         & 9 & 61,748 & 0.88 & 27.43 & 0.47 & 45,806 & 0.72 & 2.17 & 0.46 \\
\midrule
 & & \multicolumn{4}{c}{\textbf{Gurobi}} & \multicolumn{4}{c}{\textbf{Enhanced decomp heuristic + RH}} \\
\cmidrule(lr){3-6} \cmidrule(lr){7-10}
Region & \#Periods & Objective & \%Gap & \%Rel Gap & Runtime & Objective & \%Gap & \%Rel Gap & Runtime \\
\midrule
\textbf{Region 2} & 3 & \multicolumn{4}{c}{Out of memory} & 701,008 & 0.98 & / & 2.20 \\
         & 6 & \multicolumn{4}{c}{Out of memory} & 1,353,654 & 0.97 & / & 5.82 \\
         & 9 & \multicolumn{4}{c}{Out of memory} & 1,949,275 & 0.97 & / & 11.08 \\
\bottomrule
\end{tabular}
\begin{tablenotes}
\item *All runtimes are measured in hours.
\end{tablenotes}
\end{threeparttable}
}
\label{tab:heuristic_comparison}
\end{center}
\vspace{-6pt}
\end{table}

The computational results for the compact Region 1 instance are presented in Table \ref{tab:heuristic_comparison}. For all planning horizons (3, 6, and 9 periods), the direct application of Gurobi to the model yields the best-known objective values, which provide a robust benchmark. However, this level of solution quality requires significant computational effort, as Gurobi consistently terminates at the 12-hour time limit for all planning horizon lengths. It can be observed that the baseline decomposition heuristic is much faster than the Gurobi solver, solving the 9-period Region 1 instance in just 0.41 hours; However, this speed comes at a severe cost to solution quality, resulting in an objective value of \$61,743, which represents a 27.42\% relative gap compared to the Gurobi benchmark. The integration of a basic rolling horizon framework into the decomposition strategy preserves both the computational efficiency and solution quality as observed in the baseline decomposition heuristic. In contrast, our enhanced decomposition strategy with a rolling horizon and capacity-aware flow assignment significantly improves solution quality and maintains high computational efficiency. Specifically, it achieves objective values within approximately 2\% relative gap compared to the Gurobi solver across all planning horizons for Region 1, while reducing computation time to 0.11-0.46 hours. This corresponds to a speedup of roughly 25× to 100× compared to the 12-hour runtime of Gurobi. As the horizon length increases from 3 to 9 periods, the decomposition-based strategies exhibit gently near-linear scaling with the number of periods. This further demonstrates the scalability of the proposed heuristic method.

The scalability challenge is more pronounced in the large-scale Region 2 instance. As shown in Table \ref{tab:heuristic_comparison}, the Gurobi solver is out of memory even for a 3-period planning horizon, failing to produce even a feasible solution. This result highlights the practical limitations of exact methods for problems of this scale. In this context, the proposed enhanced decomposition approach with rolling horizon delivers superior performance. It successfully solves all instances for Region 2, providing solutions with low remaining optimality gaps. Moreover, its runtimes remain manageable, which scales from 2.20 hours for 3 periods to 11.08 hours for the 9-period instance. These results show that the proposed solution approach is a viable and effective tool for tackling real-world large-scale 2EVRP problems that are unsolvable by commercial solvers. Moreover, the findings underscore that decomposition is essential for memory-bounded tractability at network scale, and the proposed capacity-aware flow assignment plays a crucial role in improving solution quality.

\subsection{Network effectiveness}
\label{sec:ne}

In this section, we assess the effectiveness of different network configurations. This experiment evaluates the logistics performance of a two-tier (i.e., single-echelon) versus a three-tier (i.e., two-echelon) hyperconnected network structure under varying levels of market shares. To simulate different stages of system adoption and freight volume, we define four distinct market share scenarios, representing 0.1\%, 1\%, 10\%, and 100\% of the total freight flow within the region. The performance of each network configuration under these scenarios is measured across key logistics metrics: service capacity, operational efficiency, environmental impact (i.e., eco-impact), and detailed cost metrics. Table~\ref{tab:network_comp} summarizes the results.

The results reveal a clear and significant performance advantage for the three-tier hyperconnected network architecture as market share increases. In terms of operational efficiency, the three-tier network consistently demonstrates substantial reductions in travel time and distance, highlighting the benefit of introducing a hierarchical middle layer (i.e., local hubs) to enhance routing flexibility and reduce detours. In addition, these savings become more pronounced at higher market share levels. For example, at a 100\% market share, the three-tier network reduces total travel distance by over 0.61 million miles, a 45.0\% reduction from 1.36 million miles to 0.75 million miles, compared to the two-tier baseline. This superior routing efficiency further results in a favorable eco-impact, with CO$_2$-eq emissions reduced by nearly half. Although the two-tier network achieves a higher truck fill rate, this is primarily attributed to its limited flexibility, relying on fewer and more heavily consolidated routing patterns. In contrast, the three-tier hyperconnected network can effectively utilize the local hub tier, evidenced by the average number of LH service routes growing from 0.2 to 253 as market share scales. This structure provides greater flexibility and reduces the system's reliance on external on-demand capacity when the market share is large. As we can see from the results, the outsourcing percentage for the three-tier network becomes negligible at higher market shares.

\vspace{0.1in}
\begin{table}[h!]
\begin{center}
\caption{{Logistics performance under different market shares and network types.}}
\resizebox{\linewidth}{!}{
\begin{tabular}{llcccccccc}
\toprule
& \multicolumn{1}{l}{\textbf{\small{Market share}}} & \multicolumn{2}{c}{0.1\%} & \multicolumn{2}{c}{1\%} & \multicolumn{2}{c}{10\%} & \multicolumn{2}{c}{100\%} \\
\cmidrule(lr){3-4}\cmidrule(lr){5-6}\cmidrule(lr){7-8}\cmidrule(lr){9-10}
& \multicolumn{1}{l}{\textbf{\small{Network types}}} & 2-Tier & 3-Tier & 2-Tier & 3-Tier & 2-Tier & 3-Tier & 2-Tier & 3-Tier \\
\midrule
\textbf{\small{Service}} & \# of service routes & 53 & 67 & 570 & 624 & 5,473 & 5,535 & 51,038 & 50,372 \\
\textbf{\small{Routes}} & Avg \# of GH service routes\hspace{-0.1in} & 5 & 4 & 57 & 26 & 547 & 233 & 5,104 & 2,279 \\
& Avg \# of LH service routes\hspace{-0.1in} & / & 0.2 & / & 3 & / & 29 & / & 253 \\
\midrule
\textbf{\small{Operation}} & Travel time \midsize{(hours)} & 49 & 41 & 408 & 259 & 3,856 & 2,232 & 36,325 & 20,877 \\
\textbf{\small{Efficiency}} & Travel distances \midsize{(miles)} & 1,778 & 1,435 & 15,137 & 9,060 & 143,972 & 79,318 & 1,364,130 & 749,889 \\
& Service truck fill rate \midsize{(\%)} & 59.37 & 43.48 & 57.94 & 56.32 & 65.33 & 58.06 & 68.48 & 61.55 \\
& Outsourcing percentage \midsize{(\%)}\hspace{-0.1in} & 10.76 & 15.93 & 1.38 & 2.54 & 0.23 & 0.02 & 0.04 & 0.00 \\
\midrule
\textbf{\small{Eco-impact}}\hspace{-0.1in} 
& Energy consumption \midsize{(1,000 MJ)} 
& 11 & 9 & 90 & 54 & 854 & 470 & 8,089 & 4,446 \\
& CO$_2$-eq emissions \midsize{(metric tons)} 
& 0.8 & 0.7 & 7 & 4 & 66 & 36 & 621 & 341 \\
\midrule
\textbf{\small{Cost}} & Scheduled service cost\midsize{ (\$)} & 792 & 1,044 & 8,887 & 7,882 & 86,372 & 71,471 & 813,580 & 676,827 \\
\textbf{\small{Metrics}} & \midsize{(\%} of total cost) & (49.78) & (52.66) & (90.58) & (85.08) & (98.90) & (99.75) & (99.82) & (100.00) \\
& On-demand penalty\midsize{ (\$)}\hspace{-0.1in} & 799 & 938 & 924 & 1,382 & 964 & 180 & 1,455 & 0 \\
& \midsize{(\%} of total cost) & (50.22) & (47.34) & (9.42) & (14.92) & (1.10) & (0.25) & (0.18) & (0.00) \\
& Total cost\midsize{ (\$)} & 1,591 & 1,982 & 9,812 & 9,264 & 87,336 & 71,651 & 815,035 & 676,827 \\
\bottomrule
\end{tabular}
}
\label{tab:network_comp}
\end{center}
\vspace*{-6pt}
\end{table}

Finally, the cost structure offers additional insight into service modality tradeoffs. At a low market share of 0.1\%, the two-tier hyperconnected network shows a marginal total cost advantage, while the three-tier network becomes more cost-effective as freight volume grows. At a 10\% market share, the on-demand outsourcing cost for the three-tier system is \$180, which is 81.3\% lower than the two-tier network. At 100\% market share, the total cost savings of the three-tier network exceed \$136,000, representing a 16.8\% reduction. This cost saving is driven by more efficient routing, as the three-tier network achieves a better balance between delivery detours and freight consolidation. Moreover, at 100\% market share, the three-tier network eliminates the expensive on-demand outsourcing, which further reduces its total costs. Collectively, these results demonstrate that the three-tier hyperconnected network (i.e., two-echelon delivery system) not only shows a better performance in environmental sustainability but also achieves higher cost-effectiveness at medium to high market share levels.

\subsection{Benefits of vehicle sharing and repositioning}
\label{sec:vsr}

In this section, the analysis focuses on quantifying the operational and economic benefits derived from implementing the vehicle sharing and repositioning strategy. This experiment compares logistic performance under two strategies: one where vehicles are fixed to their initial hubs (i.e., without repositioning) and another where they can be dynamically shared and relocated across the network (i.e., with repositioning). To assess the effectiveness of the strategy under different conditions, we examine four initial vehicle distribution schemes: Case 1, uniform distribution across hubs in both echelons; Case 2, uniform distribution in the first echelon and demand-based distribution in the second; Case 3, demand-based distribution in the first echelon and uniform distribution in the second; and Case 4, demand-based distribution across both echelons. The performance of each scenario is evaluated using the same metrics as in the previous section.

\vspace{0.1in}
\begin{table}[h!]
\begin{center}
\caption{{Logistics performance under different initial vehicle distributions and repositioning options.}}
\resizebox{\linewidth}{!}{
\begin{tabular}{llcccccccc}
\toprule
& \multicolumn{1}{l}{\textbf{\small{Initial vehicle distribution}}} & \multicolumn{2}{c}{Case 1} & \multicolumn{2}{c}{Case 2} & \multicolumn{2}{c}{Case 3} & \multicolumn{2}{c}{Case 4} \\
\cmidrule(lr){3-4}\cmidrule(lr){5-6}\cmidrule(lr){7-8}\cmidrule(lr){9-10}
& \multicolumn{1}{l}{\textbf{\small{Repositioning option}}} & N & Y & N & Y & N & Y & N & Y \\
\midrule
\textbf{\small{Service}} & \# of service routes & 50,587 & 50,595 & 50,562 & 50,595 & 50,585 & 50,372 & 50,582 & 50,372 \\
\textbf{\small{Routes}} & Avg \# of GH service routes\hspace{-0.1in} & 2,287 & 2,302 & 2,287 & 2,302 & 2,287 & 2,280 & 2,287 & 2,279 \\
& Avg \# of LH service routes\hspace{-0.1in} & 254 & 253 & 254 & 253 & 254 & 253 & 254 & 253 \\
\midrule
\textbf{\small{Operation}} & Travel time \midsize{(hours)} & 24,651 & 21,830 & 24,650 & 21,828 & 24,651 & 20,878 & 24,651 & 20,877 \\
\textbf{\small{Efficiency}} & Travel distances \midsize{(miles)} & 902,311 & 789,424 & 902,313 & 789,340 & 902,318 & 749,896 & 902,312 & 749,889 \\
& Service truck fill rate \midsize{(\%)} & 55.38 & 59.85 & 55.42 & 59.86 & 55.45 & 61.52 & 55.43 & 61.55 \\
& Outsourcing percentage \midsize{(\%)}\hspace{-0.1in} & 0.22 & 0.00 & 0.22 & 0.00 & 0.22 & 0.00 & 0.22 & 0.00 \\
\midrule
\textbf{\small{Eco-impact}}\hspace{-0.1in} 
& Energy consumption \midsize{(1,000 MJ)} 
& 5,350 & 4,681 & 5,350 & 4,680 & 5,350 & 4,447 & 5,350 & 4,446 \\
& CO$_2$-eq emissions \midsize{(metric tons)} 
& 411 & 359 & 411 & 359 & 411 & 341 & 411 & 341 \\
\midrule
\textbf{\small{Cost}} & Scheduled service cost\midsize{ (\$)} & 815,127 & 712,084 & 815,127 & 712,002 & 815,149 & 676,858 & 815,142 & 676,827 \\
\textbf{\small{Metrics}} & \midsize{(\%} of total cost) & (99.12) & (100.00) & (99.12) & (100.00) & (99.12) & (100.00) & (99.12) & (100.00) \\
& On-demand penalty\midsize{ (\$)}\hspace{-0.1in} & 7,248 & 0 & 7,247 & 0 & 7,251 & 0 & 7,256 & 0 \\
&  \midsize{(\%} of total cost) & (0.88) & (0.00) & (0.88) & (0.00) & (0.88) & (0.00) & (0.88) & (0.00) \\
&  Total cost\midsize{ (\$)} & 822,375 & 712,084 & 822,374 & 712,002 & 822,400 & 676,858 & 822,398 & 676,827 \\
\bottomrule
\end{tabular}
}
\label{tab:vehicle_cases}
\end{center}
\vspace*{-6pt}
\end{table}

The results, detailed in Table \ref{tab:vehicle_cases}, demonstrate that enabling vehicle sharing and repositioning delivers consistent benefits across all initial allocation scenarios. When comparing the results with and without repositioning options, the number of scheduled routes remains nearly constant, with variations below 0.4\%. This indicates that repositioning mainly influences routing patterns and utilization rather than altering total service capacity. Moreover, the repositioning strategy reduces total travel time by 11-15\% and travel distance by 12-17\%, depending on the initial vehicle distribution. The overall truck fill rate also improves from roughly 55\% to 60\%. On the other hand, outsourcing drops to zero in every reposition-enabled case, suggesting that repositioning eliminates residual capacity gaps that previously required on-demand services. This overall efficiency gain further results in lower fuel and emissions. In Case 1, Energy consumption and CO$_2$-eq emissions are reduced by around 12.7\%. Under demand-based initial distribution (i.e., Cases 3 and 4), reductions deepen to approximately 17\%. Consequently, the repositioning strategy yields substantial economic benefits, reducing total costs by 13.4\% in Case 1 and by as much as 17.7\% in Case 4, where costs fall from \$822,398 to \$676,827. 

Furthermore, the analysis reveals that a smart initial distribution of vehicles amplifies benefits with repositioning, while the network performance is nearly identical across cases when repositioning is disabled. Specifically, aligning the initial vehicle distribution with historical demand in the first-echelon delivery (i.e., Cases 3 and 4) leads to the greatest reductions in cost, travel distance, and environmental impact. However, when comparing Case 4 with Case 3 (or Case 2 with Case 1), the results show that the initial vehicle distribution in the second echelon has far less impact than in the first echelon. Taken together, these findings highlight the effectiveness of vehicle sharing and repositioning in enhancing operational efficiency, reducing emissions, and lowering total costs, regardless of the initial vehicle distribution. The results further demonstrate the robustness of our proposed model in handling different vehicle configurations. Moreover, they provide managerial insights, suggesting that demand-aware initial fleet deployment offers an additional layer of optimization, enabling the system to achieve its potential for cost savings and sustainability.

\subsection{Environmental impacts of alternative fuel technologies}
\label{sec:ei}

This section extends the analysis to evaluate the environmental and economic impacts of deploying alternative fuel vehicle technologies within the two-echelon urban delivery network. The experiment is designed to quantify the trade-offs associated with transitioning from conventional fossil fuels to clean energy sources. We compare four distinct vehicle types based on their fuel source: (i) diesel, (ii) gasoline, (iii) battery-electric, and (iv) hydrogen (produced via natural gas reforming). The performance of each vehicle type is assessed based on three primary metrics: total cost of ownership, the well-to-wheel energy consumption, and the well-to-wheel CO$_2$-eq emissions.

\vspace{0.1in}
\begin{table}[h!]
\begin{center}
\caption{{Logistics performance under different vehicle fuel types.}}
\resizebox{0.67\linewidth}{!}{
\begin{tabular}{lcccc}
\toprule
\multicolumn{5}{c}{\textbf{\small{First-echelon urban delivery}}} \\
\midrule
\textbf{\small{Fuel Type}} 
& Diesel 
& Gasoline 
& Electricity 
& Hydrogen \\
\midrule
Cost rate\midsize{ (\$)} 
& 963,486 
& 963,486 
& 947,248 
& 1,282,844 \\
Energy consumption\midsize{ (1,000 MJ)} 
& 3,210 
& 4,169 
& 1,775 
& 2,299 \\
CO$_2$-eq emissions\midsize{ (metric tons)} 
& 246 
& 292 
& 98 
& 137 \\
\midrule
\multicolumn{5}{c}{\textbf{\small{Second-echelon urban delivery}}} \\
\midrule
\textbf{\small{Fuel Type}} 
& Diesel 
& Gasoline 
& Electricity 
& Hydrogen \\
\midrule
Cost rate\midsize{ (\$)} 
& 371,317 
& 371,317 
& 365,058 
& 494,393 \\
Energy consumption\midsize{ (1,000 MJ)} 
& 1,237 
& 1,606 
& 684 
& 886 \\
CO$_2$-eq emissions\midsize{ (metric tons)} 
& 95 
& 112 
& 38 
& 53 \\
\midrule
\multicolumn{5}{c}{\textbf{\small{\% Variation compared to Diesel truck delivery}}} \\
\midrule
\textbf{\small{Fuel Type}} 
& Diesel 
& Gasoline 
& Electricity 
& Hydrogen \\
\midrule
Cost rate\midsize{ (\%)} 
& 0.00 
& 0.00 
& -1.69 
& 33.15 \\
Energy consumption\midsize{ (\%)} 
& 0.00 
& 29.88 
& -44.69 
& -28.38 \\
CO$_2$-eq emissions\midsize{ (\%)} 
& 0.00 
& 18.49 
& -60.02 
& -44.33 \\
\bottomrule
\end{tabular}
}
\label{tab:fuel_comparison}
\end{center}
\vspace*{-6pt}
\end{table}

Table~\ref{tab:fuel_comparison} presents a detailed breakdown of cost and environmental metrics for both echelons, along with relative percentage changes benchmarked against the conventional diesel-based delivery system. The findings highlight a clear trade-off between economic costs and environmental benefits when adopting alternative fuel vehicles. Diesel and gasoline vehicles exhibit a low total cost of ownership; however, they result in poorer environmental performance, with gasoline vehicles in particular generating the highest energy consumption and emissions. Battery-electric vehicles emerge as a financially attractive option, achieving a modest cost reduction of 1.69\% relative to the diesel baseline while also delivering substantial environmental improvements. In particular, electric vehicles reduce energy consumption by 44.69\% and CO$_2$-eq emissions by 60.02\%. Hydrogen-powered vehicles, while offering great environmental benefits with a 44.33\% reduction in CO$_2$-eq emissions, are associated with considerably higher operational costs, increasing total costs by 33.15\%. Overall, these results demonstrate that the fuel choice significantly affects the environmental footprint and cost structure of urban logistics. Among the alternatives considered, battery-electric vehicles provide the most sustainable pathway toward cost savings and emission reductions, achieving near-complete decarbonization of urban delivery operations.

\subsection{Sensitivity analysis}
\label{sec:sen_analysis}

To evaluate the robustness of the proposed model and examine the influence of key cost parameters, we conduct a sensitivity analysis on the scheduled service cost and the on-demand service penalty. Seven scenarios are defined by varying the corresponding cost multipliers. Scenario 1 serves as the baseline, with both multipliers set to one (i.e., all costs at their default values). In Scenarios 2 and 3, the on-demand service cost is fixed while the scheduled service cost is doubled or halved, respectively. The default on-demand service cost is set at a very high level, evidenced by the zero outsourcing percentage in the baseline case. Therefore, we design additional scenarios by gradually reducing the on-demand service penalty to as low as one-sixteenth of its default value. The effects of these variations are assessed across multiple performance metrics, thereby providing insights into how the system’s optimal strategy adapts under different cost structures.

\vspace{0.1in}
\begin{table}[h!]
\begin{center}
\caption{{Experimental results for cost parameter sensitivity analysis.}}
\resizebox{\linewidth}{!}{ 
\begin{tabular}{llccccccc}
\toprule
& \multicolumn{1}{l}{\textbf{\small{Scheduled service cost multiplier}}}\hspace{-0.2in} & 1 & 2 & 1/2 & 1 & 1 & 1 & 1 \\
& \multicolumn{1}{l}{\textbf{\small{On-demand penalty multiplier}}}\hspace{-0.2in} & 1 & 1 & 1 & 1/2 & 1/4 & 1/8 & 1/16 \\
\midrule
\textbf{\small{Service}} & \# of service routes & 50,372 & 50,347 & 50,240 & 50,324 & 49,786 & 48,694 & 29,877 \\
\textbf{\small{Routes}} & Avg \# of GH service routes & 2,279 & 2,278 & 2,267 & 2,277 & 2,273 & 2,223 & 1,737 \\
& Avg \# of LH service routes & 253 & 253 & 253 & 253 & 248 & 243 & 115 \\
\midrule
\textbf{\small{Operation}} & Travel time \midsize{(hours)} & 20,877 & 20,873 & 20,880 & 20,873 & 20,867 & 22,345 & 62,647 \\
\textbf{\small{Efficiency}} & Travel distances \midsize{(miles)} & 749,889 & 749,748 & 749,916 & 749,740 & 749,586 & 800,341 & 2,208,263 \\
& Service truck fill rate \midsize{(\%)} & 61.55 & 61.56 & 61.25 & 61.65 & 62.14 & 63.49 & 93.97 \\
& Outsourcing percentage \midsize{(\%)} & 0.00 & 0.01 & 0.00 & 0.01 & 0.25 & 1.42 & 25.12 \\
\midrule
\textbf{\small{Eco-impact}}\hspace{-0.1in} 
& Energy consumption \midsize{(1,000 MJ)} 
& 4,446 & 4,446 & 4,447 & 4,446 & 4,445 & 4,746 & 13,094 \\
& CO$_2$-eq emissions \midsize{(metric tons)} 
& 341 & 341 & 341 & 341 & 341 & 364 & 1,005 \\
\midrule
\textbf{\small{Cost}} & Scheduled service cost\midsize{ (\$)} & 676,827 & 1,351,125 & 338,928 & 676,035 & 671,680 & 651,539 & 403,851 \\
\textbf{\small{Metrics}} & \midsize{(\%} of total cost) & (100.00) & (99.90) & (100.00) & (99.90) & (99.46) & (97.15) & (65.80) \\
& On-demand penalty\midsize{ (\$)} & 0 & 1,302 & 0 & 660 & 3,631 & 19,132 & 209,912 \\
& \midsize{(\%} of total cost) & (0.00) & (0.10) & (0.00) & (0.10) & (0.54) & (2.85) & (34.20) \\
& Total cost\midsize{ (\$)} & 676,827 & 1,352,426 & 338,928 & 676,694 & 675,311 & 670,672 & 613,762 \\
\bottomrule
\end{tabular}
}
\label{tab:service_penalty}
\end{center}
\vspace*{-6pt}
\end{table}

The results, presented in Table \ref{tab:service_penalty}, demonstrate distinct responses to variations in the two cost parameters. When the scheduled service cost is doubled or halved, the total cost of the system adjusts proportionally, increasing from \$676,827 to \$1,352,426 when doubled and decreasing to \$338,928 when halved. In contrast, the core strategy and operational metrics remain stable. For instance, the total number of service routes remains within a tight range between 50,240 and 50,372, a variation of less than 0.3\% from the baseline. Similarly, travel distances, outsourcing percentage, and environmental metrics show negligible changes. This indicates that the default on-demand penalty is greatly larger than the scheduled service cost. Therefore, the system consistently favors its internal services over on-demand outsourcing, even when the service cost is doubled.

Consequently, we reduce the on-demand penalty from its default value to one-sixteenth of the original to observe how the resulting strategies adjust under this variation. As the penalty is incrementally reduced, a clear strategic shift from scheduled services to on-demand outsourcing occurs. The outsourcing percentage rises dramatically from 0\% in the baseline case to 25.12\% when the penalty is at its lowest (i.e., 1/16). This strategic shift leads to a substantial reduction in the scheduled service's capacity, with the number of scheduled service routes dropping by over 40\% (from 50,372 to 29,877). While this leads to higher fill rates for the remaining service trucks, which increases from 61.55\% to 93.97\%. This is because low penalties drive the strategy to outsource demand whenever scheduled services are not fully utilized and exhibit low fill rates. On the other hand, we can see the operational inefficiencies due to the large amount of outsourcing. In the lowest-penalty scenario, the total travel distance nearly triples to 2,208,263 miles, and CO$_2$-eq emissions follow a similar pattern, increasing from 341 to 1,005 metric tons. This presents a trade-off: although the total cost decreases by 9.3\% when the on-demand penalty is low, this financial saving comes at the expense of a degradation in operational and environmental performance.

\section{Conclusion}
\label{sec:con}

In response to the pressing need for sustainable and efficient urban freight systems, this paper introduces and addresses a new operational problem: the two-echelon vehicle sharing and repositioning problem for containerized urban deliveries in PI-enabled hyperconnected logistics systems. We develop a novel multi-period hybrid optimization model that integrates path-based vehicle service planning with arc-based container flow routing within two-echelon delivery networks. To overcome the computational complexity inherent in large-scale urban logistics, we propose an effective solution approach that combines an echelon-region decomposition heuristic with a capacity-aware flow assignment method inside a rolling horizon framework. To the best of our knowledge, this research is the first to explicitly model and solve the integrated challenges of vehicle repositioning and containerized routing in a multi-echelon, hyperconnected urban delivery context.

We perform a comprehensive case study using real-world freight flow data from the Atlanta metropolitan area. The results validate the significant practical benefits of the proposed model and solution methodology. Specifically, our computational study demonstrates that the proposed heuristic is highly effective: for small instances, it yields near-optimal solutions with relative gaps of about 2\% in only a fraction of the runtime of commercial solvers; for large instances, it delivers solutions efficiently where commercial solvers fail to handle such scales due to memory limitations. The analysis of network structure reveals that the two-echelon (i.e., three-tier) hyperconnected delivery network consistently outperforms a single-echelon (i.e., two-tier) alternative in terms of environmental sustainability and achieves cost efficiency at medium to high market shares. At full market share, it reduces CO$_2$-eq emissions by 45.0\% and lowers total cost by 16.8\% while balancing the trade-off between freight consolidation and routing detours. Furthermore, our results underscore the benefits of implementing vehicle sharing and repositioning strategies. Compared to the cases when it is disabled, vehicle repositioning eliminates the need for on-demand services and reduces total costs by up to 17.7\%, while improving system operating efficiency and sustainability. In addition, the investigation into alternative fuels confirms that electric and hydrogen vehicles offer significant environmental benefits. Compared to traditional diesel fleets, battery-electric vehicles deliver around 60\% CO$_2$-eq reduction with a slight cost decrease, whereas hydrogen vehicles achieve 44.3\% emission reduction at a much higher cost, suggesting a practical electrification-first pathway with a hydrogen option for deep decarbonization. Finally, a sensitivity analysis demonstrates the robustness of the proposed model and provides insights highlighting that extensive use of on-demand outsourcing leads to diminished operational and environmental performance.

This research suggests several promising avenues for future investigation. First, the integration of stochastic elements, such as demand uncertainty, travel-time variability, and infrastructure availability, would extend the model to develop more resilient planning strategies. This can be achieved by incorporating the two-stage or multi-stage stochastic programs and scenario-based rolling horizons with adaptive recourse for vehicle repositioning, ensuring service reliability in dynamic urban environments. Second, it is worthwhile to explore the vehicle routing problem that incorporates both delivery and pickup operations. The consideration of backhauls and reverse logistics is expected to raise truck utilization, increase fill rates on return legs, and reduce empty miles. This extension would require more complex routing constraints to manage simultaneous capacity adjustments, but would ultimately create a more circular and efficient logistics system. Third, another valuable extension would be to co-design the two-echelon delivery problem with facility location decisions. While the present study assumes facility infrastructure as given, jointly optimizing routes and facility placement would provide a more powerful solution for long-term strategic planning. For example, one can extend the model to co-optimize charging stations or hydrogen supply sites to support phased electrification and the establishment of hydrogen corridors. Such an approach would offer insights into infrastructure investment trade-offs and guide decision-makers in transitioning toward low-carbon urban logistics systems.

\newpage
\bibliographystyle{elsarticle-harv} 
\bibliography{scibib}

\pagebreak
\appendix

\section{Complementary parameter settings}
\label{appendix_a}

Table \ref{tab:other_paras} and  \ref{tab:baseline_value} summarize the key parameters used in this study.
\begin{table}[h!]
    \begin{center}
    \caption{{Input parameter summary.}}
    \resizebox{0.73\linewidth}{!}{
    \begin{tabular}{cll}
    \toprule
    Notation & Implication & Value\\
    \midrule
    $C^{F}_{op/hr}$ & FE scheduled service cost per hour & \$40 \\
    $C^{S}_{op/hr}$ & SE scheduled service cost per hour & \$24 \\
    $C_{reposition}$ & Fixed repositioning surcharge & \$0.1 \\
    $C^P$ & On-demand service penalty per container per hour & \$50\\
    $Q_1$ & Capacity of a FE vehicle & 10 containers \\
    $Q_2$ & Capacity of a SE vehicle & 4 containers \\
    $V$ & Average delivery vehcile speed & 33 miles/hour\\
    $H$ & length of rolling horizon windows & 3 periods \\
    $N$ & Number of operational periods per day & 3 \\
    \bottomrule
    \label{tab:other_paras}
    \end{tabular}
    }
    \end{center}
    \vspace{-0.25in}
\end{table}

\vspace{-0.05in}
\begin{table}[h!]
\begin{center}
\caption{{Well-to-wheels economic and environmental parameters for different vehicle fuel types \citep{Basma2023TCOClass8, GREET2025Updates}.}}
\resizebox{0.76\linewidth}{!}{
\begin{tabular}{lccccc}
\toprule
\textbf{\small{Fuel Type}} & Diesel & Gasoline & Electricity & Hydrogen \\
\midrule
Per-mile cost rate\midsize{ (\$ per mile)} & 1.78 & 1.78 & 1.75 & 2.37  \\
Per-mile energy consumption\midsize{ (MJ per mile)} & 5.93 & 7.70 & 3.28 & 4.25 \\
Per-mile CO$_2$-eq emissions\midsize{ (kg per mile)} & 0.46 & 0.54 & 0.18 & 0.25 \\
\bottomrule
\end{tabular}
}
\label{tab:baseline_value}
\end{center}
\vspace*{-16pt}
\end{table}

\section{Route set generation and pre-processing}
\label{appendix_b}

The proposed optimization model relies on a predefined set of service routes, denoted as $R$ for the first echelon and $M$ for the second echelon. Given the substantial number of potential routes in a complex network, we employ a two-stage methodology to construct a comprehensive yet computationally manageable set of candidate routes. First, a large number of potential routes is generated by enumeration based on the considered route patterns. Second, this initial set is reduced through a multi-step pre-processing procedure designed to filter for efficiency while ensuring network coverage.

Based on the operational structure of the two-echelon delivery system, we have two distinct sets of service routes: 1) first-echelon routes, which facilitate the movement of goods between gateway hubs and transfer them to the downstream local hub layer, and 2) second-echelon routes, which handle distribution from local hubs to the final access hubs in delivery zones. To manage the problem complexity, the number of intermediate stops in any single route is limited. Specifically, we consider four FE route patterns: i) \(g\!-\!g\), a repositioning point-to-point route between two gateway hubs; ii) \(g\!-\!g-\!g\), a non-repositioning point-to-point route between two gateway hubs; iii) \(g\!-\!l\!-\!g\), a route visiting one local hub between two gateways; and iv) \(g\!-\!l\!-\!l\!-\!g\), a route visiting two distinct local hubs between its origin and destination gateways. Regarding the SE route, we consider the following two route patterns: \(l\!-\!a\!-\!l\), a route connecting two local hubs via a single access hub; and \(l\!-\!a\!-\!a\!-\!l\), a route connecting two local hubs via two distinct access hubs. It is noted that FE route patterns iii) and iv) and SE route patterns i) and ii) include both repositioning and non-repositioning routes, depending on whether the starting and ending hubs are the same. Then, given these route patterns, we are able to generate potential route sets by enumerating the network arcs. 

However, this raw enumeration generates a large number of inefficient paths with long detours that are dominated by alternative paths. To eliminate such candidates and improve the computational tractability, a multi-step pre-processing methodology is implemented. Specifically, we first apply a thresholding rule at the segment level to filter the inefficient routes. To establish these filtering thresholds, we analyze the distance distribution of each arc type (e.g., \(g\text{--}g\), \(l\text{--}l\), \(a\text{--}a\), \(g\text{--}l\), \(l\text{--}a\)). We then define a maximum permissible distance threshold for each arc type based on a specific percentile of its distribution. For instance, the threshold for short-haul \(a\text{--}a\) segments is set at a low percentile (e.g., 10th) to select dense and local connections. Correspondingly, each candidate route in the initial set is evaluated. If any consecutive arc in its path exceeds the distance threshold for its specific type, the route is removed from the candidate set. In practice, we find that this segment-wise rule is effective in removing inefficient paths while preserving sensible routes.

Nevertheless, the above filtering method can potentially remove all viable routes for certain hubs, isolating them from the network. To address this issue and ensure model integrity, we apply a coverage restoration procedure. This restoration guarantees that each active hub in the network is served by a minimum set of required route types operating within its designated area. For access hubs $a$, we ensure the existence of at least one second–echelon tour \(l\!-\!a\!-\!l\) whose local hub lies in $a$’s assigned local cell. For local hubs $l$, we ensure at least one customer-serving tour in $\{l\!-\!a\!-\!l,\; l\!-\!a\!-\!a\!-\!l\}$ and at least one first–echelon route \(g\!-\!l\!-\!g\) with the gateway hub within the same urban area. For gateway hubs $g$, we ensure the presence of a satellite-serving path in $\{g\!-\!l\!-\!g,\; g\!-\!l\!-\!l\!-\!g\}$ whose local hub is in $g$’s urban area. The algorithm iterates through each hub and verifies if these route types exist for it in the reduced set. If a required route type is missing for a given hub, the procedure first searches the original, unfiltered route set for the most efficient candidate (i.e., the one with the shortest travel time) that satisfies the requirement and adds it back into the reduced set. In the rare case where no such route exists in the original set, a non-repositioning round-trip route is generated, connecting the hub to the nearest required hub type within the same area. Overall, the proposed route generation, segment-aware pruning, and coverage restoration prove to be effective in our empirical study. It effectively reduces the solution space without sacrificing too much solution quality.

\end{document}